# RATES OF CONVERGENCE FOR THE POSTERIOR DISTRIBUTIONS OF MIXTURES OF BETAS AND ADAPTIVE NONPARAMETRIC ESTIMATION OF THE DENSITY

BY JUDITH ROUSSEAU[1]

*Université Paris Dauphine and CREST*

In this paper, we investigate the asymptotic properties of nonparametric Bayesian mixtures of Betas for estimating a smooth density on $[0,1]$. We consider a parametrization of Beta distributions in terms of mean and scale parameters and construct a mixture of these Betas in the mean parameter, while putting a prior on this scaling parameter. We prove that such Bayesian nonparametric models have good frequentist asymptotic properties. We determine the posterior rate of concentration around the true density and prove that it is the minimax rate of concentration when the true density belongs to a Hölder class with regularity $\beta$, for all positive $\beta$, leading to a minimax adaptive estimating procedure of the density. We also believe that the approximating results obtained on these mixtures of Beta densities can be of interest in a frequentist framework.

**1. Introduction.** In this paper, we study the asymptotic behaviour of posterior components. There is a vast literature on mixture models because of their rich structure which allows for different uses; for instance, they are well known to be adapted to the modelling of heterogeneous populations as is used, for example, in cluster analysis (for a good review on mixture models see [10] or [11] for various aspects of Bayesian mixture models). They are also useful in nonparametric density estimation, in particular, they can be considered to capture small variations around a specific parametric model, as typically occurs in robust estimation or in a goodness of fit test of a parametric family or of a specific distribution (see, e.g., [12, 13]). The approach considered here is density estimation, but it has applications in many other aspects of mixture models, such has clustering, classification and goodness

Received September 2008; revised February 2009.
[1]Supported in part by the ANR-SP Bayes grant.

of fit testing, since in all of these cases, understanding the behaviour of the posterior distribution is crucial. Nonparametric prior distributions based on mixture models are often considered in practice and Dirichlet mixture priors are particularly popular. Dirichlet mixtures have been introduced by [2, 9] and have been widely used ever since, but their asymptotic properties are not well known apart from a few cases such as Gaussian mixtures, triangular mixtures and Bernstein polynomials. The papers [4, 5] and [15] study the concentration rate of the posterior distribution under Dirichlet mixtures of Gaussian priors, and Ghosal [3] considers the Bernstein polynomial's case, that is, the mixture of Beta distribution with fixed parameters. The paper [13] considers mixtures of triangular distributions, with a prior on the mixing distribution which is not necessarily a Dirichlet process. In all of those cases, the authors mainly consider the concentration rate of the posterior around the true density, when the latter have some known regularity conditions, or when it is a continuous mixture.

Posterior distributions associated with Bernstein polynomials are known to be suboptimal in terms of minimax rates of convergence when the true density is Hölder. An improvement is obtained in [8] based on a modification of Bernstein polynomials leading to the minimax rate of convergence in the classes of Hölder densities with regularity $\beta$ when $\beta \leq 1$. In this paper, we consider another class of mixtures of Beta models which is richer and, therefore, allows for better asymptotic results.

Beta densities are often represented as

$$(1.1) \qquad g(x|a,b) = \frac{x^{a-1}(1-x)^{b-1}}{B(a,b)}, \qquad B(a,b) = \frac{\Gamma(a)\Gamma(b)}{\Gamma(a+b)}.$$

Here we consider a different parametrization of the Beta distribution writing $a = \alpha/(1-\varepsilon)$ and $b = \alpha/\varepsilon$, so that $\varepsilon \in (0,1)$ is the mean of the Beta distribution, and $\alpha > 0$ is a scale parameter. To approximate smooth densities on $[0,1]$, we consider a location mixture of Beta densities in the form,

$$(1.2) \quad g_{\alpha,P}(x) = \sum_{j=1}^{k} p_j g_{\alpha,\varepsilon_j}(x), \qquad g_{\alpha,\varepsilon_j}(x) = g(x|\alpha/(1-\varepsilon_j), \alpha/\varepsilon_j),$$

where the mixing density is given by

$$(1.3) \qquad P(\varepsilon) = \sum_{j=1}^{k} p_j \delta_{\varepsilon_j}(\varepsilon).$$

The parameters of this mixture model are then $k \in N^*$, and for each $k$, $(\alpha, p_1, \ldots, p_k, \varepsilon_1, \ldots, \varepsilon_k)$. The prior probability on the set of densities can, therefore, be expressed as

$$d\pi(f) = p(k)\pi_k(\varepsilon_1, \ldots, \varepsilon_k, p_1, \ldots, p_k|\alpha) \, d\pi_{k,\alpha}(\alpha), \qquad \text{if } f = g_{\alpha,P},$$

or $d\pi(f) = d\pi(P|\alpha) \, d\pi_2(\alpha)$, in the case of a Dirichlet mixture.

Determining the concentration rate of the posterior distribution around the true density corresponds to determining a sequence $\tau_n$ converging to 0 such that if

$$(1.4) \qquad B_{\tau_n} = \{f \in \mathcal{F}, d(f, f_0) < \tau_n\}$$

for some distance or pseudo-distance $d(\cdot, \cdot)$ on the set of densities, and if $X^n = (X_1, \ldots, X_n)$ where the $X_i$'s are independent and identically distributed from a distribution having a density $f_0$ with respect to Lebesgue measure, then

$$(1.5) \qquad P^\pi[B_{\tau_n}|X^n] \to 1, \qquad \text{in probability.}$$

The difficulty with mixture models comes from the fact that it is often quite hard to obtain precise approximating properties for these models. The papers [14, 18] give general descriptions of the Kullback–Leibler support of priors based on mixture models. These results are key in obtaining the consistency of the posterior distribution, but cannot be applied to obtain rates of concentration. In these papers, they use the Kernel structure of mixture models. Among such mixture models, location-scale kernels are widely considered. mixtures of Betas are not location-scale kernels. However, when $\alpha$ gets large, $g_{\alpha,\varepsilon}$ concentrates around $\varepsilon$, so that locally, these Beta densities behave like Gaussian densities. This behavior is described in Section 3. Using these ideas, we study the approximation of a density $f$ by a continuous mixture in the form

$$(1.6) \qquad g_{\alpha,f}(x) = \int_0^1 f(\varepsilon) g_{\alpha,\varepsilon}(x) \, d\varepsilon,$$

where $f$ is a probability density on $(0,1)$. When $\alpha$ becomes large, $g_{\alpha,\varepsilon}(x)$ behaves locally like a location scale kernel so that $g_{\alpha,f}$ becomes close to $f$. Similarly to the Gaussian case, this approximation is *good* only if $f$ has a regularity less than 2. However, by shifting slightly the mixing density, it is possible to improve the approximation so that continuous mixtures of Betas are *good* approximations of any smooth density (see Section 3.1). As in the case of Gaussian mixtures (see [4, 15]), we approximate the continuous mixture by a discrete mixture. In [5], the authors derive a posterior rate of concentration of the posterior distribution around the true density when the true density is twice continuously differentiable. In particular, they obtain the minimax rate $n^{-2/5}$, up to a $\log n$ term under the $L_1$ risk.

In this paper, we show that the minimax rate can be obtained (up to a $\log n$ term) for any $\beta > 0$ by choosing carefully the rate at which $\alpha$ increases with $n$ and considering a prior on $\alpha$ leads to an adaptive minimax rate of concentration of the posterior. This result has much theoretical and practical interest.

1.1. *Notation.* Throughout the paper, $X_1, \ldots, X_n$ are independent and identically distributed as $P_0$, having density $f_0$, with respect to Lebesgue measure. We assume that $X_i \in [0,1]$. We consider the following three distances (or pseudo-distances) on the set of densities on $[0,1]$: the $L_1$ distance: $\|f - g\|_1 = \int_0^1 |f(x) - g(x)| \, dx$, the Kullback–Leibler divergence: $\mathrm{KL}(f, g) = \int_0^1 f(x) \log(f(x)/g(x)) \, dx$, for any densities $f, g$ on $[0,1]$ and for any $k > 1$ $V_k(f, g) = \int_0^1 f(x) \times |\log(f(x)/g(x))|^k \, dx$. We also denote by $\|g\|_\infty$ the supremum norm of the function $g$.

$\mathcal{H}(L, \beta)$ denotes the class of Hölder functions with regularity function $\beta$: let $r$ be the largest integer smaller than $\beta$, and denote by $f^{(r)}$ its $r$th derivative.

$$\mathcal{H}(L, \beta) = \{f : [0,1] \to \mathbb{R}; |f^{(r)}(x) - f^{(r)}(y)| \leq L|x - y|^{\beta - r}\}.$$

We denote by $\mathcal{S}_k$ the simplex, $\mathcal{S}_k = \{y \in [0,1]^k; \sum_{i=1}^k y_i = 1\}$.

We denote by $P^\pi[\cdot|X^n]$, the posterior distribution given the observations $X^n = (X_1, \ldots, X_n)$, and $E^\pi[\cdot|X^n]$, the expectation with respect to this posterior distribution. Similarly $E_0^n$ and $P_0^n$ represent the expectation and the probability with respect to the true density $f_0^{\otimes n}$ and $E_f^n$ and $P_f^n$ the expectation and probability with respect to the distribution $f^{\otimes n}$.

1.2. *Assumptions.* Throughout the paper, we assume that the true density $f_0$ is positive on the open interval $(0, 1)$ and satisfies:

ASSUMPTION $\mathbf{A}_0$. If $f_0 \in \mathcal{H}(\beta, L)$, there exist integers $0 \leq k_0, k_1 < \beta$ such that

$$f^{(k_0)}(0) > 0, \qquad f^{(k_1)}(1) < 0;$$

$k_0$ and $k_1$ denote the first integers such that the corresponding derivatives calculated at 0 and 1, respectively, are nonzero.

This assumption is quite mild and ensures that $f_0(x)$ does not go too quickly to 0 when $x$ goes to 0 or 1 so that we can control the Kullback–Leibler divergence between $f_0$ and mixtures of Betas.

1.3. *Organization of the paper.* The paper is organized as follows. In Section 2, we give the two main theorems on the concentration rates of the posterior distributions under specific types of priors. In Section 3, we present some results describing the approximating properties of mixtures of Betas. We believe that these results are interesting outside the Bayesian framework, since they could also be applied to obtain convergence rates for maximum likelihood estimators. This section is divided into two parts. First we describe how continuous mixtures can approach smooth densities (Section 3.1), then we approach continuous mixtures by discrete mixtures (Section 3.2). Finally, Section 4 is dedicated to the proofs of Theorems 2.1 and 2.2.

**2. Posterior concentration rates.** In this section, we give the two main results on the concentration rates of the posterior distribution around the true density. We first consider the case of a varying number of components, which we call the *adaptive prior* and then we consider a Dirichlet mixture also leading to an adaptive rate of concentration on a more restrictive class of densities. In both cases, a diffuse prior on $\alpha$ is considered. Finally, a nonadaptive rate is obtained by considering a deterministic sequence $\alpha_n$ increasing to infinity. We consider a concentration rate in terms of the $L_1$ distance, however, the results can be applied to the Hellinger distance as well. We first describe the *adpative prior*.

*Adaptive prior*: let $f = g_{\alpha,P}$ and $P = \sum_{i=1}^{k} p_i \delta_{\varepsilon_i}$ the mixing distribution then

$$d\pi(f) = p(k)\, d\pi_{k,2}(p_1,\ldots,p_k) \prod_{j=1}^{k} \pi_e(\varepsilon_j)\pi_\alpha(\alpha)\, d\alpha\, d\varepsilon_1 \cdots d\varepsilon_k.$$

For all $k > 0$, $d\pi_{k,2}$ has a positive density $\pi_{k,2}$ with respect to Lebesgue measure on the simplex $\mathcal{S}_k$, which is bounded from below by a term in the form $c_1^k$. Conditionally on $k$, the $\varepsilon_j$'s, $j = 1,\ldots,k$, are independent and identically distributed with density $\pi_e$ which satisfies

$$a_1 \varepsilon^T (1-\varepsilon)^T \geq \pi_e(\varepsilon) \geq a_2 \varepsilon^T (1-\varepsilon)^T \qquad \forall \varepsilon \in (0,1),$$

for some $a_1, a_2 > 0$, and $T \geq 1$. We consider the following conditions on the prior $\pi_\alpha$: $\pi_\alpha$ is bounded and for all $b_1 > 0$, there exist $c_1, c_2, c_3, A > 0$ such that for all $u$ large enough,

$$\pi_\alpha(c_1 u < \alpha < c_2 u) \geq C e^{-b_1 u^{1/2}},$$
$$\pi_\alpha(c_3 u < \alpha) \leq C e^{-b_1 u^{1/2}},$$
$$\pi_\alpha(\alpha < e^{-uA}) \leq C e^{-b_1 u}.$$

Let $L(k)$ be either equal to 1 for all $k$ or $L(k) = \log(k)$. The distribution on $k$ satisfies the following condition: there exist $a_1, a_2 > 0$ such that for all $K$ large enough,

$$e^{-a_1 \mathrm{KL}(K)} \leq p[k = K] \leq e^{-a_2 \mathrm{KL}(K)}.$$

Note that if $\sqrt{\alpha}$ follows a Gamma distribution with parameters $(a, b)$ with $a \geq 1$, then the conditions on $\pi_\alpha$ are satisfied. We have the following theorem:

THEOREM 2.1. *Consider an* adaptive *prior, as described above, then the posterior distribution satisfies, for all $\beta > 0$ and $f_0 \in \mathcal{H}(\beta, L)$ satisfying Assumption* $\mathbf{A}_0$,

$$P^\pi[B^c_{\tau_n}|X^n] = o_P(1)$$

*with*

$$\tau_n = \tau_0 n^{-\beta/(2\beta+1)}(\log n)^{5\beta/(4\beta+2)}, \qquad \text{if } L(k) = \log(k),$$
$$\tau_n = \tau_0 n^{-\beta/(2\beta+1)}(\log n)^{5\beta/(4\beta+2)+1/2}, \qquad \text{if } L(k) = 1.$$

The prior does not depend on $\beta$ so that the procedure is adaptive and optimal up to a $\log n$ term, since for each $\beta > 0$ the rate $n^{-\beta/(2\beta+1)}$ is the minimax rate of convergence in the class $\mathcal{H}(\beta, L)$.

Dirichlet mixtures form an alternative to the above prior, which is often considered in practice, since they lead to efficient algorithms and have interesting properties for classification models, for instance. We now present the asymptotic concentration rate of the posterior based on the following Dirichlet mixtures of Beta densities.

*Dirichlet prior*: the mixing distribution $P$ follows a Dirichlet process $\mathcal{D}(\nu)$ associated with a finite measure whose density with respect to Lebesgue measure is denoted $\nu$ and is positive on the open interval $(0, 1)$. Assume also that $\nu$ is bounded and satisfies

$$\nu(\xi) \geq \nu_0 \xi^{T_1}(1 - \xi)^{T_1}.$$

The prior on $\alpha$, $\pi_\alpha$ has support $[n^t, +\infty)$, for some $0 < t < 1$ and satisfies, for all $b_1 > 0$, there exist $c_1, c_2, c_3, C > 0$ such that for all $\alpha_n$ satisfying $\alpha_n n^{-t} \to +\infty$

$$\pi_\alpha(c_1 \alpha_n < \alpha < c_2 \alpha_n) \geq C e^{-b_1 \sqrt{\alpha_n}},$$
$$\pi_\alpha(c_3 \alpha_n < \alpha) \leq C e^{-b_1 \sqrt{\alpha_n}}.$$

Note that if $\sqrt{\alpha} \stackrel{d}{=} n^t + \Gamma(a, b)$, with $a, b > 0$, then the above condition is satisfied.

THEOREM 2.2. *Consider a* Dirichlet *prior then the posterior distribution satisfies: for all $f_0 \in \mathcal{H}(\beta, L)$ with $\beta > 0$, and satisfying Assumption* $\mathbf{A}_0$,

$$P^\pi[B^c_{\tau_n}|X^n] = o_P(1)$$

*with*

$$\tau_n = \tau_0 n^{-\beta/(2\beta+1)}(\log n)^{5\beta/(2\beta+1)}, \qquad \text{if } \beta \leq 1/t - 1/2,$$
$$\tau_n = \tau_0 n^{-1/2+t/4}(\log n)^{(6\beta+1/2)/(2\beta+1)}, \qquad \text{if } \beta > 1/t - 1/2.$$

Hence the *Dirichlet prior* implies a minimax adaptive rate of concentration on the densities with regularity $\beta < 1/t - 1/2$. By choosing $t$ small, this class of functions is quite large, with small loss in the rates of convergence.

We could have considered $\alpha = \alpha_n$ deterministic and increasing with $n$, which would have implied the following nonadaptive posterior rate, depending on $\alpha_n$.

COROLLARY 2.1. *Consider a prior belonging either to the class of* adaptive priors *or to the class of* Dirichlet prior, *as described above, apart from the fact that* $\alpha = \alpha_n = o(n)$ *is deterministic. Then if* $f_0 \in \mathcal{H}(\beta, L)$, *and satisfies Assumption* $\mathbf{A}_0$,

$$P^\pi[B^c_{\tau_n}|X^n] = o_P(1)$$

$$\text{with } \tau_n = \tau_0(\log \alpha_n)[\alpha_n^{-\beta/2} \vee (\sqrt{\alpha_n \log \alpha_n}/n)^{1/2}].$$

In particular, if $\alpha_n = n^{2/(2\beta+1)}(\log n)^{-3/(2\beta+1)}$, we obtain the minimax rate (up to a $\log n$) term $\tau_n = \tau_0 n^{-\beta/(2\beta+1)}(\log n)^{5\beta/(4\beta+2)}$. Note that deterministic sequences $\alpha_n$ lead to nonadaptive concentration rates.

These results imply that for any $\beta > 0$, the optimal rate, in the minimax sense, is obtained. Hence the above mixtures of Betas form a richer class of models than the Bernstein polynomials or the mixtures of triangular distributions who lead, at best, to the minimax rates for $\beta \leq 2$. It is to be noted, however, that Bernstein polynomials and mixtures of triangular densities have other interesting properties and are particularly easy to simulate.

Corollary 2.1 sheds light on the impact of $\alpha_n$ as a scale parameter. It can thus be compared to the scale parameter $\sigma_n$ which appears in Dirichlet mixtures of Gaussian distributions. Interestingly, van der Vaart and van Zanten [16, 17] also study the impact of scaling factors in nonparametric priors constructed as scaled Gaussian processes, and as in our case, considering a random scaling factor allows for adaptive, minimax concentration rates.

In Section 3, we see that the key factor leading to such a rate is the possibility of approximating any $f_0 \in \mathcal{H}(L, \beta)$ by a continuous mixture in the form $g_{\alpha_n, f}$ with an error of order $\alpha_n^{-\beta}$, for some density $f$ close to $f_0$ but not necessarily equal to $f_0$. An interesting feature leading to this approximating property is that $g_{\alpha_n, \varepsilon}$ acts locally as a Gaussian kernel around $\varepsilon$. However, the interest in the Bayesian procedure, compared to a classical frequentist kernel nonparametric method, comes from the fact that we do not necessarily need to approach $f_0$ by $g_{\alpha_n, f_0}$, which would have constrained us to $\beta \leq 2$. Indeed, if necessary, we can consider a slight modification $f$ of $f_0$ such that $g_{\alpha_n, f}$ approximates $f_0$ with an error of order $\alpha_n^{-\beta}$ for all $\beta$. This is described in the following section.

**3. Approximation of a smooth density by continuous and discrete mixtures.** A Beta mixture, as defined by (1.6) behaves locally like a Gaussian mixture, however, its behaviour seems to be richer since the variance adapts to the value of $x$ (see Lemma 3.1). In this section, we obtain a way to approximate any Hölder density $f$ by a sequence of continuous and discrete mixtures. We begin with approximating the density by a sequence of continuous mixtures, and then we approximate the continuous mixtures by discrete mixtures.

3.1. *Continuous mixtures.* We consider a continuous mixture $g_{\alpha,f}$ as defined in (1.6). This mixture is based on the parametrization of a beta density in terms of mean $\varepsilon$ and scale $\alpha$. The idea in this section is that when $\alpha$ becomes large, the above mixture converges to $f$, if $f$ is continuous. We first give a result where the approximation is controlled in terms of the supremum norm, which has an intrinsic interest. We also give a bound on the approximation error for Kullback–Leibler-types of divergence, which is the required result to control the posterior concentration rate.

THEOREM 3.1. *Assume that $f_0 \in \mathcal{H}(\beta, L)$ and satisfies Assumption $\mathbf{A}_0$, with $\beta > 0$. Then there exists a probability density $f_1$ such that*

$$f_1(x) = f_0(x)\left(1 + \sum_{j=2}^{\lceil\beta\rceil-1} \frac{w_j(x)}{\alpha^{j/2}}\right), \qquad \text{if } \beta > 2;$$

$$f_1(x) = f_0(x), \qquad \text{if } \beta \leq 2,$$

*where the $w_j$'s are combinations of polynomial functions of $x$ and of terms in the form $f_0^{(l)}(x) x^l (1-x)^l / f_0(x)$, $l \leq j$, and*

$$\|g_{\alpha,f_1} - f_0\|_\infty \leq C\alpha^{-\beta/2}, \tag{3.1}$$

*and for all $p > 0$,*

$$\text{KL}(f_0, g_{\alpha,f_1}) \leq C\alpha^{-\beta}, \qquad \int f_0 \left|\log\left(\frac{f_0}{g_{\alpha,f_1}}\right)\right|^p \leq C\alpha^{-\beta}. \tag{3.2}$$

Note that if we do not allow $f_1$ to be different from $f_0$, we do not achieve the rate $\alpha^{-\beta}$ to be true for values of $\beta$ greater than 2. We believe that the trick of allowing $f_1$ to be different from $f_0$ could be used in a more general context of Bayesian mixture distributions (or Bayesian kernel approaches as defined in [18]), inducing a greater flexibility of Bayesian kernel methods with respect to frequentist kernel methods.

A Beta density with parameters $(\alpha/\varepsilon, \alpha/(1-\varepsilon))$ can be expressed as

$$g_{\alpha,\varepsilon}(x) = x^{\alpha/(1-\varepsilon)-1}(1-x)^{\alpha/\varepsilon-1}\frac{\Gamma(\alpha/(\varepsilon(1-\varepsilon)))}{\Gamma(\alpha/\varepsilon)\Gamma(\alpha/(1-\varepsilon))}.$$

From this, we have the following three approximations that will be used throughout the proofs of Theorems 2.1, 2.2, 3.1 and 3.2. Let

$$K(\varepsilon, x) = \varepsilon \log(\varepsilon/x) + (1-\varepsilon)\log((1-\varepsilon)/(1-x)), \tag{3.3}$$

this is the Kullback–Leibler divergence between the Bernoulli $\varepsilon$ and the Bernoulli $x$ distributions. Then:

LEMMA 3.1.

$$g_{\alpha,\varepsilon}(x) = \frac{\sqrt{\alpha}}{\sqrt{2\pi x(1-x)}} e^{-\alpha K(\varepsilon,x)/(\varepsilon(1-\varepsilon))}$$
$$\times \left[1 + \sum_{j=1}^{k} \frac{b_j(\varepsilon)}{\alpha^j} + O(\alpha^{-(k+1)})\right] \tag{3.4}$$

for any $k > 0$ and $\alpha$ large enough where the $b_j(\varepsilon)$ are polynomial functions. For all $k > 0, k_1 \geq 3$, we also have,

$$g_{\alpha,\varepsilon}(x) = \frac{\sqrt{\alpha}}{\sqrt{2\pi x(1-x)}}$$
$$\times \exp\left\{-\frac{\alpha(x-\varepsilon)^2}{2x^2(1-x)^2}\right.$$
$$\left.\times \left[1 + \frac{(x-\varepsilon)}{x(1-x)}\left(C(x) + Q_{k_1}\left(\frac{x-\varepsilon}{x(1-x)}\right)\right)\right] + R_1\right\} \tag{3.5}$$
$$\times \left[1 + \sum_{j=1}^{k} \frac{b_j(\varepsilon)}{\alpha^j} + O(\alpha^{-(k+1)})\right],$$

where $R_1 \leq \alpha C|x-\varepsilon|^{k_1-2}(x_\varepsilon(1-x_\varepsilon))^{-k_1+2}$,

$$Q_{k_1}\left(\frac{x-\varepsilon}{x(1-x)}\right) = \sum_{l=0}^{k_1-3} \frac{C_l(x)(x-\varepsilon)^l}{(x(1-x))^l}$$

and the functions $C(x), C_l(x), l \leq k_1$, are polynomial where $x_\varepsilon \in (x, \varepsilon)$ and $C$ is a positive constant. Moreover when $\alpha|x-\varepsilon|^3 \leq C_0 x^3(1-x)^3$ for any positive constant $C_0$, if $k_2 \geq 0$, and if $k_1 \geq 3 \vee 3k_2$, there exists $C_1 > 0$ such that

$$g_{\alpha,\varepsilon}(x) = \frac{\sqrt{\alpha} e^{-\alpha(x-\varepsilon)^2/(2x^2(1-x)^2)}}{\sqrt{2\pi x(1-x)}}$$
$$\times \left(\sum_{j=0}^{k_2} \frac{\alpha^j (x-\varepsilon)^{3j}}{j!(x(1-x))^{3j}}\left[C(x) + Q_{k_1}\left(\frac{x-\varepsilon}{x(1-x)}\right)\right]^j + R\right) \tag{3.6}$$

$$\times \left[ 1 + \sum_{j=1}^{k} \frac{b_j(\varepsilon)}{\alpha^j} + O(\alpha^{-(k+1)}) \right],$$

where $|R| \leq C_1 \alpha^{k_2+1} |x - \varepsilon|^{3(k_2+1)} (x_\varepsilon(1-x_\varepsilon))^{-3(k_2+1)}$.

Note that the term $O(\alpha^{-(k+1)})$ appearing in (3.4), (3.5) and (3.6) is uniform in $x$ and $\varepsilon$.

PROOF OF LEMMA 3.1. The proof of (3.4) follows from the expression of the Beta densities in the form,

$$g_{\alpha,\varepsilon}(x) = \frac{\Gamma(\alpha/(\varepsilon(1-\varepsilon)))\varepsilon^{\alpha/(1-\varepsilon)}(1-\varepsilon)^{\alpha/\varepsilon}}{\Gamma(\alpha/\varepsilon)\Gamma(\alpha/(1-\varepsilon))} \frac{e^{-\alpha K(\varepsilon,x)/(\varepsilon(1-\varepsilon))}}{x(1-x)}$$

and from a Taylor expansion of $\Gamma(y)$ for $y$ close to infinity where we obtain that

$$\frac{\Gamma(\alpha/(\varepsilon(1-\varepsilon)))}{\Gamma(\alpha/\varepsilon)\Gamma(\alpha/(1-\varepsilon))}$$
$$= \frac{\sqrt{\alpha}}{\sqrt{2\pi}} \exp\left(-\alpha \left[\frac{\log(\varepsilon)}{1-\varepsilon} + \frac{\log(1-\varepsilon)}{\varepsilon}\right]\right) \left(1 + \sum_{j=1}^{\infty} b_j \frac{\varepsilon^j(1-\varepsilon)^j}{\alpha^j}\right)$$
$$\times \left(1 + \sum_{j=1}^{\infty} b_j \frac{\varepsilon^j}{\alpha^j}\right)^{-1} \left(1 + \sum_{j=1}^{\infty} b_j \frac{(1-\varepsilon)^j}{\alpha^j}\right)^{-1},$$

where the $b_j$'s are the coefficient appearing in the expansion of the Gamma function near infinity (see, e.g., [1]). Putting the three remaining terms together results in: for all $k > 0$,

$$\left(1 + \sum_{j=1}^{\infty} b_j \frac{\varepsilon^j(1-\varepsilon)^j}{\alpha^j}\right) \left(1 + \sum_{j=1}^{\infty} b_j \frac{\varepsilon^j}{\alpha^j}\right)^{-1} \left(1 + \sum_{j=1}^{\infty} b_j \frac{(1-\varepsilon)^j}{\alpha^j}\right)^{-1}$$
$$= 1 + \sum_{j=1}^{k} \frac{b_j(\varepsilon)}{\alpha^j} + O(\alpha^{-(k+1)}),$$

where the $b_j(\varepsilon)$'s are polynomial functions with degree less than $2j$. This implies (3.4). To obtain (3.5) we make a Taylor expansion of (3.4) as a function of $\varepsilon$ around $x$.

$$\frac{K(\varepsilon,x)}{\varepsilon(1-\varepsilon)} = \frac{(\varepsilon-x)^2}{2x^2(1-x)^2} + \sum_{j=3}^{k_1} C_j(x) \frac{(x-\varepsilon)^j}{x^j(1-x)^j} + R_1,$$

where $R_1 \leq R|x - \varepsilon|^{k_1+1}/(x_\varepsilon(1-x_\varepsilon))^{k_1+1}$ for some $x_\varepsilon \in (x,\varepsilon)$, leading to (3.5). A Taylor expansion of $e^y$ around 0 combined with the above approximation of $y$ leads to (3.6). □

To prove (3.1), we control the difference between the uniform density on $[0,1]$ and the corresponding Beta mixture $g_\alpha = \int_0^1 g_{\alpha,\varepsilon}\,d\varepsilon$. This is given in the following lemma.

LEMMA 3.2. *For all $\alpha > 0$ large enough, for all $k_2 \geq 1$ and $k_1 \geq 3(k_2 - 1)$ define*

$$I(x) = \sum_{j=1}^{k_2} \frac{C(x)^j \mu_{3j}}{\alpha^{j/2}} + \sum_{l=2}^{k_2 k_1} \frac{B_l(x)}{\alpha^{l/2}} \mu_j = E[\mathcal{N}(0,1)^j],$$

*then*

$$\left\| g_\alpha(x) - 1 - \frac{I(x)}{\alpha} \right\|_\infty \leq C\alpha^{-(k_2+1)/2} (\log \alpha)^{3(k_2+1)/2},$$

*where the $B_l(x)$'s are polynomial functions of $x$.*

The proof of Lemma 3.2 is given in Appendix A. We now prove Theorem 3.1.

PROOF OF THEOREM 3.1. Throughout the proof, $C$ denotes a generic positive constant. Let $f \in \mathcal{H}(\beta, L)$ and denote $r = \lfloor \beta \rfloor$. Then $\forall \varepsilon \in (0,1)$,

$$(3.7) \qquad \left| f(\varepsilon) - \sum_{j=0}^{r} \frac{f^{(j)}(x)}{j!}(\varepsilon - x)^j \right| \leq L|x - \varepsilon|^\beta.$$

The construction of $f_1$ is iterative. Let $\delta_x = \delta_0 x(1-x)\sqrt{\log \alpha/\alpha}$. We bound

$$\int_0^1 |x-\varepsilon|^\beta g_{\alpha,\varepsilon}(x)\,d\varepsilon \leq \left| \int_0^{x-\delta_x} g_{\alpha,\varepsilon}(x)\,d\varepsilon + \int_{x+\delta_x}^1 g_{\alpha,\varepsilon}(x)\,d\varepsilon \right|$$
$$+ \int_{x-\delta_x}^{x+\delta_x} |x-\varepsilon|^\beta g_{\alpha,\varepsilon}(x)\,d\varepsilon.$$

Equation (A.6) implies that for all $H > 0$, if $\delta_0$ is large enough, the first term of the right-hand side of the above inequality is $O(\alpha^{-H})$. We treat the second term using the same calculations as in the case of $I_3$ in Appendix A, so that for all $k > 0$,

$$\int_{x-\delta_x}^{x+\delta_x} |x-\varepsilon|^\beta g_{\alpha,\varepsilon}(x)\,d\varepsilon$$
$$\leq C\alpha^{-\beta/2} x^\beta (1-x)^\beta E[|\mathcal{N}(0,1)|^\beta] + O(\alpha^{-k/2}).$$

Therefore,

$$\int_0^1 |x-\varepsilon|^\beta g_{\alpha,\varepsilon}(x)\,d\varepsilon = O(\alpha^{-\beta/2} x^\beta (1-x)^\beta) + O(\alpha^{-H}) \quad \forall H > 0,$$

uniformly in $x$. Then for all $H > 0$,

$$[g_{\alpha,f} - f](x) = \sum_{j=1}^r \frac{f^{(j)}(x)}{j!} \int_0^1 (\varepsilon - x)^j g_{\alpha,\varepsilon}(x)\,d\varepsilon + f(x)(g_\alpha(x) - 1)$$

$$+ O(\alpha^{-\beta/2} x^\beta (1-x)^\beta) + O(\alpha^{-H})$$

$$= \sum_{j=1}^r \frac{f^{(j)}(x)}{j!} \int_0^1 (\varepsilon - x)^j g_{\alpha,\varepsilon}(x)\,d\varepsilon + f(x) \frac{I(x)}{\alpha}$$

$$+ O(\alpha^{-\beta/2} x^\beta (1 - x^\beta) + \alpha^{-H}),$$

uniformly in $x$, for all $H > 0$. Using the same calculations as in the computation of $I_3$ in the proof Lemma 3.2, we obtain for all $j \geq 1$, to the order $O(\alpha^{-(k+j+1)/2} x^j (1-x)^j + \alpha^{-H})$

$$\int_0^1 (\varepsilon - x)^j g_{\alpha,\varepsilon}(x)\,d\varepsilon$$

$$= \frac{\sqrt{\alpha}}{\sqrt{2\pi} x(1-x)} \int_{x-\delta_x}^{x+\delta_x} e^{-\alpha(x-\varepsilon)^2/(2x^2(1-x)^2)}$$

$$\times \left( (x-\varepsilon)^j + \sum_{l=1}^k \frac{\alpha^l (x-\varepsilon)^{3l+j}}{j!(x(1-x))^{3l}} \left[ C(x) + Q_{k_1}\left(\frac{x-\varepsilon}{x(1-x)}\right) \right]^l \right) d\varepsilon$$

$$= \mu_j \alpha^{-j/2} x^j (1-x)^j + \sum_{l=1}^k \frac{D_l(x) x^j (1-x)^j}{\alpha^{(j+l)/2}},$$

so that we can write,

$$\int_0^1 (\varepsilon - x)^j g_{\alpha,\varepsilon}(x)\,d\varepsilon$$

$$= \frac{x^j(1-x)^j}{\alpha^{j/2}} \mu_{j,\alpha}(x) + O(\alpha^{-(k+j+1)/2} x^j (1-x)^j + \alpha^{-H}),$$

where $\mu_{j,\alpha}(x)$ is a polynomial function of $x$ with the leading term being equal to $\mu_j$. We can thus write, to the order $O(\alpha^{-\beta/2} x^\beta (1 - x^\beta) + \alpha^{-H})$

$$(3.8) \quad [g_{\alpha,f} - f](x) = \sum_{j=1}^r \frac{f^{(j)}(x) x^j (1-x)^j \mu_{j,\alpha}(x)}{j! \alpha^{j/2}} + f(x) \frac{I(x)}{\alpha}.$$

Hence if $\beta \leq 2$, since $\mu_1 = 0$,

$$|g_{\alpha,f} - f|(x) \leq \frac{\|I\|_\infty f(x)}{\alpha} + O(\alpha^{-\beta/2} x^\beta (1-x^\beta)) + O(\alpha^{-H})$$
(3.9)
$$= O(\alpha^{-\beta/2})$$

as soon as $H > \beta/2$, leading to (3.1) with $f_1 = f$. If $\beta > 2$, we construct a probability density $f_1$ satisfying

$$(g_{\alpha,f_1} - f)(x) = O(\alpha^{-\beta/2} x^\beta (1-x)^\beta) + O(\alpha^{-H}).$$

Equation (3.8) implies that $f_1$ needs satisfy, to the order $O(\alpha^{-H})$,

$$\sum_{j=1}^{r} \frac{f_1^{(j)}(x) x^j (1-x)^j \mu_{j,\alpha}(x)}{j! \alpha^{j/2}} + f_1(x)\left(1 + \frac{I(x)}{\alpha}\right)$$

$$= f(x) + O(\alpha^{-\beta/2} x^\beta (1-x^\beta)).$$

To prove that such a probability density exists we construct it iteratively. Let $2 < \beta \leq 3$, then set

$$h_1(x) = f(x)\left(1 - \frac{I(x)}{\alpha}\right) - \frac{x(1-x) f'(x) C(x) \mu_4}{\alpha} - \frac{x^2 (1-x)^2 f''(x) \mu_2}{2\alpha}.$$

Note that if $f \in \mathcal{H}(L, \beta)$, then $\inf f > 0$ implies $h_1 > 0$ for $\alpha$ large enough and if $f(0) = 0$ [$f(1) = 0$] when $x$ is close to 0 (resp. 1), if

$$\liminf_x \frac{f(x)}{x^j (1-x)^j |f^{(j)}(x)|} > 0, \qquad j = 1, 2,$$

$h_1 \geq 0$ for $\alpha$ large enough on $[0,1]$. Assumption $\mathbf{A}_0$ implies the above relation between $f$ and $f^{(j)}$ since

$$h_1(x) = \frac{x^{k_0} f^{(k_0)}(\bar{x}_1)}{k_0!}\left(1 - \frac{I(x)}{\alpha}\right) - \frac{x^{k_0}(1-x) f^{(k_0)}(\bar{x}_2) C(x) \mu_4}{\alpha(k_0 - 1)!}$$
(3.10)
$$- \frac{x^{k_0}(1-x)^2 f^{(k_0)}(\bar{x}_3) \mu_2}{2\alpha(k_0 - 2)!}$$

with $\bar{x}_1, \bar{x}_2, \bar{x}_3 \in (0, x)$. Since $f^{(k_0)}(0) > 0$, $h_1(x)$ is equivalent to $f(x)$ for $\alpha$ large enough and $x$ close to zero, and $h_1(x) > 0$ for all $x \in (0,1)$. Let $c_1 = \int_0^1 h_1(x)\,dx$. Since $\int_0^1 [g_{\alpha,f} - f](x)\,dx = 0$,

$$c_1 = 1 + O(\alpha^{-(3/2 \wedge \beta/2)})$$

and we can divide $h_1$ by its normalizing constant and obtain the same result as before, so that $h_1$ can be chosen to be a probability density on $[0, 1]$.

From this we obtain when $\beta > 2$,

$$(g_{\alpha,h_1} - f)(x) = \int_0^1 \left( \sum_{j=1}^{r-2} \frac{h_1^{(j)}(x)}{j!}(\varepsilon - x)^j \right) g_{\alpha,\varepsilon}(x) \, d\varepsilon + h_1(x) \frac{I(x)}{\alpha}$$

$$+ \sum_{j=r-1}^{r} \left( f(x) - \frac{I(x)}{\alpha} \right)^{(j)} \frac{\int_0^1 (\varepsilon - x)^j g_{\alpha,\varepsilon}(x) \, d\varepsilon}{j!}$$

$$+ O(\alpha^{-\beta/2} x^\beta (1-x)^\beta)$$

$$= \frac{w(x) f(x)}{\alpha^2} + O(\alpha^{-2 \wedge \beta/2} x^\beta (1-x)^\beta) + O(\alpha^{-H}) \quad \forall H > 0,$$

where $w(x)$ is a combination of polynomial functions of $x$ and of functions in the form $x^j(1-x)^j f^{(j)}(x)$, with $j < 3$, if $\beta \leq 4$. If $\beta \leq 4$, then we set $f_1 = h_1$ (renormalized), else we reiterate. We thus obtain that if $r_\beta$ is the largest integer (strictly) smaller than $\beta/2$,

$$f_1(x) = f(x) \left( 1 + \sum_{j=1}^{\lfloor r_\beta \rfloor} \frac{w_j(x)}{\alpha^j} \right),$$

where $w_j(x)$ is a combination of polynomial functions and of terms in the form $f^{(l)}(x) x^l (1-x)^l / f(x)$, $l \leq 2j$. Assumption $\mathbf{A}_0$ implies that $f_1$ can be chosen to be a density when $\alpha$ is large enough and satisfies

$$\|g_{\alpha,f_1} - f\|_\infty \leq C \alpha^{-\beta/2},$$

which implies (3.1).

If $f$ is strictly positive on $[0,1]$, then (3.2) follows directly from (3.1). We now consider the case where $f(0) = 0$ [the case $f(1) = 0$ is treated similarly]. Under the Assumption $\mathbf{A}_0$, the previous calculations lead to

$$(g_{\alpha,f_1} - f)(x) = O(f(x) \alpha^{-\beta/2}) + O(\alpha^{-H}) \quad \forall H > 0.$$

Note also that for $\alpha$ large enough, $f_1$ is increasing between 0 and $\delta$ for some positive constant $\alpha > 0$, so that if $x$ is small enough,

(3.11)
$$g_{\alpha,f_1} \geq \frac{f_1(x) \sqrt{\alpha}}{2\sqrt{2\pi} x(1-x)} \int_x^{x+\delta_x} e^{-\alpha(x-\varepsilon)^2/(2x^2(1-x)^2)} \, d\varepsilon$$

$$\geq \frac{f_1(x)}{4},$$

so that $g_{\alpha,f_1} \geq f/8$ on $[0,1]$. Therefore, since $f(x) = f^{(k_0)}(0) x^{k_0}/k_0! + o(x^{k_0})$ when $x$ is close to 0, let $H > \beta$ and $c = c_0 \alpha^{-H/k_0}$; for some constant $c_0$ large

enough, we have

$$\mathrm{KL}(f, g_{\alpha, f_1}) \leq \log 2 \int_0^c f(x)\, dx + \alpha^{-\beta} \int_c^1 f(x)\, dx$$

$$+ \int_c^1 f(x) \left| \log\left(1 - \frac{\alpha^{-H}}{f(x)}\right) \right| dx$$

$$\leq C(\alpha^{-H(k_0+1)/k_0} + \alpha^{-\beta} + \alpha^{-H}) = O(\alpha^{-\beta}).$$

Similarly, for all $p > 0$, if $c_p = c_0 \alpha^{-H/(pk_0)}$,

$$\int f(x) |\log(f(x)/g_{\alpha, f_1}(x))|^p\, dx$$

$$\leq (\log 2)^p \int_0^{c_p} f(x)\, dx + \alpha^{-p\beta} \int_{c_p}^1 f(x)\, dx$$

$$+ \int_{c_p}^1 f(x) \left| \log\left(1 - \frac{\alpha^{-H}}{f(x)}\right) \right|^p dx$$

$$\leq C(\alpha^{-2H(k_0+1)/(pk_0)} + \alpha^{-p\beta} + \alpha^{-H})$$

$$= O(\alpha^{-\beta}),$$

if $H \geq p\beta$. This achieves the proof of Theorem 3.1. $\square$

In the following section, we consider the approximation of continuous mixtures by discrete mixtures in a way similar to [4].

3.2. *Discrete mixtures.* Let $P$ be a probability on $[0,1]$ with cumulative distribution function denoted by $P(x)$ for all $x \in [0,1]$. We consider a mixture of Betas similar to before but with general probability distribution $P$ on $[0,1]$,

$$g_{\alpha, P}(x) = \int_0^1 g_{\alpha, \varepsilon}(x)\, dP(\varepsilon).$$

Let $f$ be a probability density with respect to Lebesgue measure on $[0,1]$. In this section, we study the approximation of $g_{\alpha, f}$ by $g_{\alpha, P}$ where $P$ is a discrete measure with finite support.

The approximation of discrete mixtures by continuous ones is studied in different contexts of location scale mixtures. See, for instance, [4] or [7], Chapter 3, for a general result. Beta mixtures are not location scale mixtures, however, as discussed in the previous section; when $\alpha$ is large they behave locally like location scale mixtures. In this section, we use this property to approximate continuous mixtures with finite mixtures having a reasonably small number of points in their support.

THEOREM 3.2. *Let $f$ be a probability density on $[0,1]$, $f(x) > 0$ for all $0 < x < 1$, and such that there exists $k_1, k_0 \in \mathbb{N}$ satisfying $f(x) \sim x^{k_0} c_0$, if $x = o(1)$ and $f(1-x) \sim (1-x)^{k_1} c_1$, if $1-x = o(1)$. Then there exists a discrete probability distribution $P$ having at most $N = N_0 \sqrt{\alpha} (\log \alpha)^{3/2}$ points in its support, such that for all $p \geq 1$, for all $H > 0$ (depending on $M_0$), for $\alpha$ large enough,*

$$(3.12) \qquad \int_0^1 g_{\alpha,f} \left| \log\left(\frac{g_{\alpha,f}}{g_{\alpha,P}}\right) \right|^p (x) \, dx \leq C \alpha^{-H}.$$

*We can choose the distribution $P$ such that there exists $A > 0$ with $p_j > \alpha^{-A}$ for all $j \leq N$.*

We use this inequality to obtain the following result on the true density $f_0$.

COROLLARY 3.1. *Let $f_0 \in \mathcal{H}(L, \beta)$, $\beta > 0$, be a probability density on $[0,1]$, satisfying $f_0(x) > 0$ for all $0 < x < 1$, and such that there exist $k_1, k_0 \in \mathbb{N}$ satisfying $|f^{(k_0)}(0)| > 0$ and $|f^{(k_1)}(1)| > 0$, $k_0, k_1 < \beta$. Then for all $p > 1$, there exists a discrete probability distribution $P$ having at most $N = N_0 \sqrt{\alpha} (\log \alpha)^{3/2}$ in its support, with $N_0$ large enough such that*

$$(3.13) \qquad \mathrm{KL}(f_0, g_{\alpha,P}) \leq C\alpha^{-\beta}, \qquad V_p(f_0, g_{\alpha,P}) \leq C\alpha^{-\beta}.$$

PROOF. From Theorem 3.1 there exists $f_1$ positive with $f_1 = f_0(1 + O(\alpha^{-1}))$ and

$$\mathrm{KL}(f_0, g_{\alpha,f_1}) \leq C\alpha^{-\beta}, \qquad g_{\alpha,f_1} \geq f_0/8.$$

This implies that

$$\mathrm{KL}(f_0, g_{\alpha,P}) \leq \mathrm{KL}(f_0, g_{\alpha,f_1}) + \left| \int f_0(x) \log(g_{\alpha,f_1}/g_{\alpha,P})(x) \, dx \right|$$

$$\leq C\alpha^{-\beta} + 8 \int g_{\alpha,f_1}(x) \left| \log\left(\frac{g_{\alpha,f_1}}{g_{\alpha,P}}\right) \right| (x) \, dx = O(\alpha^{-\beta}).$$

The same calculations apply to $\int_0^1 f_0(x) |\log(f_0(x)/g_{\alpha,P}(x))|^p \, dx \leq C\alpha^{-\beta}$, which achieves the proof of Corollary 3.1. □

PROOF OF THEOREM 3.2. The proof follows the same line as in [4], except that we have to control the approximations in places where the Gaussian approximation to Betas cannot be applied. Throughout this proof, $C$ denotes a generic positive constant. We first bound the difference between both mixtures at all $x$. By symmetry, we can consider $x \in [0, 1/2]$. Consider the following approximation of the exponential: for all $s \geq 0$ and all $z > 0$,

$$
(3.14) \qquad \left| e^{-z} - \sum_{j=0}^{s} \frac{(-1)^j z^j}{j!} \right| \leq \frac{z^{s+1}}{(s+1)!}.
$$

Equation (3.5) implies that for all $k > 1, k_1 \geq 3$, there exist polynomial functions of $x$, $D_l(x)$, $l \leq k_1$, and polynomial functions of $\varepsilon$, $b_j(\varepsilon)$, $j \leq k$, such that if $|x - \varepsilon| \leq M\delta_0\sqrt{\log \alpha}x(1-x)/\sqrt{\alpha}$, and setting

$$
0 \leq z = \frac{(x-\varepsilon)^2 \alpha}{2x^2(1-x)^2}\left(1 + \sum_{l=1}^{k_1-2} \frac{D_l(x)(x-\varepsilon)^l}{x^l(1-x)^l}\right) \leq CM^2 \log(\alpha),
$$

$$
g_{\alpha,\varepsilon}(x) = \frac{\sqrt{\alpha}e^{R_{k_1}}}{\sqrt{2\pi}x(1-x)}e^{-z}\left(1 + \sum_{j=1}^{k} \frac{b_j(\varepsilon)}{\alpha^j} + O(\alpha^{-(k+1)})\right),
$$

where $|R_{k_1}| \leq \alpha C \alpha^{-k_1/2+1/2}(\log \alpha)^{k_1/2}$. Consider $\varepsilon_0 = \alpha^{-t_0}$, for some positive constant $t_0$ and $\varepsilon_j = \varepsilon_0(1 + M\sqrt{\log \alpha}/\sqrt{\alpha})^j$, $j = 1, \ldots, J$, with

$$
J = \left\lfloor \frac{t_0 \log(\alpha) + 2\log(\log(\alpha))}{\log(1 + M\sqrt{\log \alpha}/\sqrt{\alpha})} \right\rfloor + 1 = O(\sqrt{\alpha}\sqrt{\log \alpha}).
$$

Define $dF_j$ and $dP_j$ the renormalized probabilities $dF$ and $dP$ restricted to $[\varepsilon_j, \varepsilon_{j+1})$ set $H > 0$. Consider $k_1 - 1 > 2H$ and $k \geq H - 1/2$ and $x \in [\varepsilon_{j-1}, \varepsilon_{j+2}]$, $j \geq 2$, using (3.14) together with the above approximation of $g_{\alpha,\varepsilon}$, we consider the moment matching approach of [4] (Lemma A.1) so that we can construct a discrete probability $dP_j$ with at most $N = 2kk_1s + 1$ supporting points such that for all $l \leq 2sk_1, l' \leq k$,

$$
\int \varepsilon^l b_{l'}(\varepsilon) d(F_j - P_j)(\varepsilon) = 0,
$$

leading to

$$
\left| \int_{\varepsilon_j}^{\varepsilon_{j+1}} g_{\alpha,\varepsilon}(x)[dF_j - dP_j](\varepsilon) \right|
$$

$$
(3.15) \qquad \leq \frac{C}{x(1-x)}\left[\frac{\sqrt{\alpha}C^{s+1}M^{2(s+1)}\log \alpha^{s+1}}{(s+1)!} + \alpha^{-H}\right]
$$

$$
= \frac{O(\alpha^{-H})}{x(1-x)},
$$

if $s = s_0 \log \alpha$ with $s_0 \geq C^2 M^4 + 1$ and $s_0 \log s_0 > H$. Moreover, for all $x \leq \varepsilon_{j-1}$, using (3.4) and the fact that

$$
\alpha K(\varepsilon, x) \geq \alpha K(\varepsilon, \varepsilon_j) \geq \frac{M^2(\log \alpha)\varepsilon(1-\varepsilon)}{3},
$$

when $\varepsilon_{j+1} > \varepsilon > \varepsilon_j$, we obtain

$$g_{\alpha,\varepsilon}(x) \leq \frac{C}{x(1-x)} e^{-cM^2 \log \alpha}$$

for some positive constant $c > 0$. A similar argument implies that if $x > \varepsilon_{j+2}$

$$g_{\alpha,\varepsilon}(x) \leq \frac{C}{x(1-x)} e^{-cM^2 \log \alpha}$$

for some positive constant $c > 0$. Hence, by constructing $P$ in the form, if $\varepsilon_{J+2} = 1 - \varepsilon_0$

$$dP(\varepsilon) = \sum_{j=0}^{J} (F(\varepsilon_{j+1}) - F(\varepsilon_j)) \, dP_j(\varepsilon) + F(\varepsilon_0) \delta_{(\varepsilon_0)} + (1 - F(\varepsilon_{J+2})) \delta_{(\varepsilon_{J+2})},$$

we finally obtain for all $x$,

$$(3.16) \qquad \left| \int_0^1 g_{\alpha,\varepsilon}(x) [dF - dP](\varepsilon) \right| \leq \frac{C\alpha^{-H}}{x(1-x)},$$

where $P$ has at most $N_\alpha = N_0 (\log \alpha)^{3/2} \sqrt{\alpha}$, for some $N_0 > 0$ related to $H$. We now consider $x \leq \varepsilon_0 (1 - M\sqrt{\log \alpha / \alpha})$. We use the approximation (3.4).

$$g_{\alpha,\varepsilon_0}(x) = \frac{C\sqrt{\alpha}}{x(1-x)} e^{-\alpha K(\varepsilon_0,x)/(\varepsilon_0(1-\varepsilon_0))} (1 + O(\alpha^{-1})).$$

Since, when $x \leq \varepsilon_0$,

$$\frac{K(\varepsilon_0, x)}{\varepsilon_0(1-\varepsilon_0)} \leq (1-\varepsilon_0)^{-1} (\log(\varepsilon_0/x)),$$

we obtain

$$g_{\alpha,P}(x) \geq e^{-\alpha \log(\varepsilon_0/x)} \frac{C\sqrt{\alpha} F(\varepsilon_0)}{x(1-x)}$$

and using the above inequalities on $g_{\alpha,\varepsilon}(x)$ for $x < \varepsilon_{j-1}$, we have

$$g_{\alpha,P}(x) \leq C\alpha^{-H}/x(1-x),$$

where $H$ depends on $M$, so that

$$|\log(g_{\alpha,P}(x))| \leq C\alpha |\log(x)|.$$

Since $g_{\alpha,f}$ is bounded (as a consequence of the fact that $g_{\alpha,f} - f$ is uniformly bounded whenever $f$ is continuous), and since $u|\log(u)|^p$ goes to zero when $u$ goes to zero,

$$(3.17) \qquad \int_0^{\varepsilon_0} g_{\alpha,f}(x) \left| \log \left( \frac{g_{\alpha,f}(x)}{g_{\alpha,P}(x)} \right) \right|^p dx \leq C\alpha^{-t_0} + C\alpha^{-t_0} (\log \alpha)^p$$
$$= O(\alpha^{-t_0} (\log \alpha)^p).$$

Note also that if $\alpha$ is large enough,
$$g_{\alpha,f}(x) \geq f(x)/4,$$
so that $g_{\alpha,f}(x) \geq cx^{k_0}(1-x)^{k_1}$ for $x$ close to 0 and for all $x \in (\varepsilon_0, 1-\varepsilon_0)$, for all $H > 0$
$$\frac{|g_{\alpha,f}(x) - g_{\alpha,P}(x)|}{g_{\alpha,f}(x)} \leq C \frac{\alpha^{-H}}{x^{k_0+1}(1-x)^{k_1+1}} \leq C\alpha^{-H+t_0(1+k_0 \vee k_1)}.$$

So that if $H > t_0(1+k_0 \vee k_1) + B/p$, with $B > 0$,
$$\int_0^1 g_{\alpha,f}(x) \left|\log\left(\frac{g_{\alpha,f}(x)}{g_{\alpha,P}(x)}\right)\right|^p dx \leq C\alpha^{-t_0}(\log \alpha)^p + C\alpha^{-B} = O(\alpha^{-B})$$

as soon as $t_0 > B$. Moreover, we can assume that there exists a fixed $A$ such that for all $j$, $p_j > \alpha^{-A} = v$. Indeed, let $I_v = \{j; p_j \leq v\}$, then consider for $j \notin I_v$, $\tilde{p}_j = cp_j$ and for $j \in I_v$, $\tilde{p}_j = cv$ where $c$ is defined by $\sum_{j=1}^J \tilde{p}_j = 1$. This implies in particular that
$$|c-1| \leq vJ \leq J_0 \alpha^{-A+1/2}(\log \alpha)^{3/2}.$$

Let $\tilde{P} = \sum_{j=0}^J \tilde{p}_j \delta_{\varepsilon_j}(\varepsilon)$ then $g_{\alpha,\tilde{P}} \geq cg_{\alpha,P}$ and if $A - 1/2 > B$,
$$\mathrm{KL}(g_{\alpha,f}, g_{\alpha,\tilde{P}}) \leq C\alpha^{-B} + |\log c| \leq C'\alpha^{-B}.$$

Also,
$$\int |g_{\alpha,\tilde{P}} - g_{\alpha,P}| \leq \alpha^{-A+1/2}(\log \alpha)^{3/2},$$

hence, if $A$ is large enough, inequality (3.16) is satisfied with $\tilde{P}$ instead of $P$. Since $p_0 = F_1(\varepsilon_0) \geq F_0(\varepsilon_0)/4$ and $F_0(\varepsilon_0) \geq \alpha^{-t_0 k_0} C$, by choosing $A > t_0 k_0$, we obtain that $0 \notin I_v$ and
$$g_{\alpha,\tilde{P}}(x) \geq g_{\alpha,\varepsilon_0}(x) F(\varepsilon_0) \qquad \forall x < \varepsilon_0,$$

so that (3.17) is satisfied with $\tilde{P}$ instead of $P$, which leads to: for all $B > 0$ there exists a distribution $\tilde{P}$, having less than $N_0 \sqrt{\alpha}(\log \alpha)^{3/2}$ points in its support, satisfying: $\tilde{p}_j \geq \alpha^{-A}$ for some $A > 0$ and all $j$ and such that
$$\int g_{\alpha,f} \left|\log\left(\frac{g_{\alpha,f}}{g_{\alpha,\tilde{P}}}\right)\right|^p (x) \, dx = O(\alpha^{-B}),$$

which achieves th proof of Theorem 3.2. □

Note, however, that $A$ depends on $B$ and so does $N_0$. Note also that this result could be used to obtain a rate of concentration of the posterior distribution around the true density when the latter is a continuous mixture.

In the following sections, we give the proofs of Theorems 2.1 and 2.2.

**4. Proofs of Theorems 2.1 and 2.2.** To prove these theorems we use Theorem 4 of [6]. In particular, let $p \geq 2$, and following their notation define

$$B^*(f_0, \tau, p) = \{f; \mathrm{KL}(f_0, f) \leq \tau^2; V_p(f_0, f) \leq \tau^p\}.$$

We also denote $\mathcal{J}_n(\tau) = N(\tau, \mathcal{F}_n, \|\cdot\|_1)$, the $L_1$ metric entropy on the set $\mathcal{F}_n$, that is, the logarithm of the minimal number of balls with radii $\tau$ needed to cover $\mathcal{F}_n$ where $\mathcal{F}_n$ is a set of densities that will be defined in each of the proofs. The proofs consist in obtaining a lower bound on $\pi(B^*(f_0, \tau_n, p))$ and an upper bound on $\mathcal{J}_n(\tau_n)$ when $f_0$ belongs to $\mathcal{H}(\beta, L)$.

4.1. *Proof of Theorem 2.1.* Assume that $f_0 \in \mathcal{H}(\beta, L)$ and let $\tau_n = \alpha_n^{-\beta/2}$, with $\alpha_n$ an increasing sequence to infinity. We first bound from below $\pi(B^*(f_0, \tau_n, p))$. Let $\alpha \in (c_1 \alpha_n, c_2 \alpha_n)$, $0 < c_1 < c_2$, using Corollary 3.1 there exists a probability distribution with $N_n = N_0 \sqrt{\alpha} (\log \alpha)^{3/2}$ supporting points such that

$$\mathrm{KL}(f_0, g_{\alpha, P}) \leq C \alpha^{-\beta}, \qquad V_p(f_0, g_{\alpha, P}) \leq C \alpha^{-\beta},$$

with $P$ of the form,

$$P(\varepsilon) = \sum_{j=1}^{k_n} p_j \delta_{\varepsilon_j}(\varepsilon),$$

$\varepsilon_j \in (\alpha^{-\beta}(\log \alpha)^{-\beta-1}, 1 - \alpha^{-\beta}(\log \alpha)^{-\beta-1})$ and $p_j > \alpha^{-A}$ for all $j = 1, \ldots, N_n$ and some fixed positive constant $A$. Set $\varepsilon_0 = \alpha^{-\beta}(\log \alpha)^{-\beta-1}$, then $\varepsilon_1 > \varepsilon_0$. Consider $dP'(\varepsilon) = \sum_{j=1}^{k_n} p'_j \delta_{\varepsilon'_j}(\varepsilon)$ with $|\varepsilon'_j - \varepsilon_j| \leq a\alpha^{-\gamma_1} \varepsilon_j (1 - \varepsilon_j)$ and $|p_j - p'_j| \leq a\alpha^{-\gamma_1 + 1/2} p_j$, for some positive constant $\gamma_1 > 1/2$. Note that this implies that $|p'_j - p_j| \leq 2a\alpha^{-\gamma_1 + 1/2} p'_j$. Then

$$(4.1) \qquad \mathrm{KL}(f_0, g_{\alpha, P'}) \leq C\alpha^{-\beta} + \int_0^1 f_0(x) \log\left[\frac{g_{\alpha, P}(x)}{g_{\alpha, P'}(x)}\right] dx.$$

For the purpose of symmetry, we work on $x \leq 1/2$. Let $M_n = M\sqrt{\log \alpha}/\sqrt{\alpha}$, when $|x - \varepsilon_j| \leq M_n \varepsilon_j (1 - \varepsilon_j)$, then Lemma B.1 implies that

$$\left|\frac{g_{\alpha, \varepsilon_j}}{g_{\alpha, \varepsilon'_j}}(x) - 1\right| = O(\alpha^{-(\gamma_1 - 1/2)} \sqrt{\log \alpha}),$$

by choosing $k_2 > 2\gamma_1 - 1$ and $k_3 > \gamma_1 - 1/2$. Set $\beta_0 > 0$, then for all $x > e^{-\beta_0 \alpha_n}$ and all $j'$ such that $|x - \varepsilon'_j| > M_n \varepsilon_j (1 - \varepsilon_j)$; since $\varepsilon_j (1 - \varepsilon_j) \geq \alpha^{-t_0}/2$ with $t_0 > \beta$, Lemma B.1 implies that if $\gamma_1 > t_0 + \beta + 2$

$$\left|\frac{g_{\alpha, \varepsilon_j}}{g_{\alpha, \varepsilon'_j}}(x) - 1\right| \leq C\beta_0 \alpha^{-\gamma_1 + 2} \varepsilon_0^{-1} = O(\alpha^{-\beta}).$$

This implies that if $x \in (e^{-\beta_0 \alpha}, 1 - e^{-\beta_0 \alpha})$,

$$\frac{g_{\alpha,P}(x)}{g_{\alpha,P'}(x)} = 1 + \frac{\sum_{j=1}^{k_n}(p_j - p'_j)g_{\alpha,\varepsilon_j}}{\sum_{j=1}^{k_n} p'_j g_{\alpha,\varepsilon'_j}} + \frac{\sum_{j=1}^{k_n} p'_j(g_{\alpha,\varepsilon_j} - g_{\alpha,\varepsilon'_j})}{\sum_{j=1}^{k_n} p'_j g_{\alpha,\varepsilon'_j}}$$

(4.2)

$$= 1 + O(\alpha^{-\gamma_1+2}).$$

Now let $x < e^{-\beta_0 \alpha}$, then $|x - \varepsilon_j| \geq \varepsilon_j(1 - \varepsilon_j)/2$ for all $j = 0, \ldots, N_n$, and there exists $c > 0$ independent of $\beta_0$ such that

$$g_{\alpha,P}(x) \leq e^{-c\alpha} \frac{\sqrt{\alpha}}{x(1-x)}, \qquad g_{\alpha,P'}(x) \leq e^{-c\alpha} \frac{\sqrt{\alpha}}{x(1-x)}.$$

Note also that

$$g_{\alpha,\varepsilon}(x) \geq C \frac{\sqrt{\alpha}}{x(1-x)} e^{-\alpha K(\varepsilon,x)/(\varepsilon(1-\varepsilon))},$$

where

$$\frac{\alpha K(\varepsilon, x)}{\varepsilon(1-\varepsilon)} = \alpha \left( \frac{1}{1-\varepsilon} \log(\varepsilon/x) + \frac{1}{\varepsilon} \log(1-\varepsilon) + \frac{x}{\varepsilon} + o(x/\varepsilon) \right)$$

$$\leq \alpha \left( \frac{1}{1-\varepsilon} \log(\varepsilon/x) + \frac{1}{\varepsilon} \log(1-\varepsilon) + \frac{x}{\varepsilon} \right) + o(1).$$

Consider the function

$$h(\varepsilon) = \frac{1}{1-\varepsilon} \log(\varepsilon/x) + \frac{1}{\varepsilon} \log(1-\varepsilon) + \frac{x}{\varepsilon},$$

since $x < |\log(1-\varepsilon)|$ for all $\varepsilon \in (\varepsilon_0, 1-\varepsilon_0)$ $h$ is increasing, and for all $\varepsilon < 1/2$, $h(\varepsilon) \leq 2|\log(x)| + O(1)$. This leads to

(4.3) $$g_{\alpha,P}(x) \geq CP([0, 1/2]) \frac{\sqrt{\alpha}}{x(1-x)} e^{2\alpha \log(x)}.$$

The same inequality holds for $g_{\alpha,P'}$, which implies that

$$\int_0^{e^{-\beta_0 \alpha}} f_0(x) \left| \log\left(\frac{g_{\alpha,P}(x)}{g_{\alpha,P'}(x)}\right) \right|^p dx \leq C\alpha^{p+1} e^{-\beta_0 \alpha} \qquad \forall p \geq 1.$$

The same kind of inequalities are obtained for $x > 1 - e^{-\beta_0 \alpha}$. Finally we obtain

(4.4) $$\int_0^1 f_0(x) \left| \log\left(\frac{g_{\alpha,P}(x)}{g_{\alpha,P'}(x)}\right) \right|^p dx = O(\alpha^{-\beta}).$$

Note that if $|p'_j - p_j| \leq \alpha^{-\beta-A}$, then $|p'_j - p_j| \leq \alpha^{-\beta} p_j$, so we need only determine a lower bound on the prior probability of the following set under the *adaptive prior*: set $\beta_0 < 1/2$

$$S_n = \{\underline{p}' \in \mathcal{S}_{N_n}; |p'_j - p_j| \leq \alpha_n^{-\beta-A}, j \leq N_n\}$$
$$\times \{|\varepsilon_j - \varepsilon'_j| \leq \alpha_n^{-2\beta-1}\varepsilon_j(1-\varepsilon_j), j \leq N_n\}.$$

The prior probability of $S_{n,1} = \{\underline{p}' \in \mathcal{S}_{N_n}; |p'_j - p_j| \leq \alpha_n^{-\beta-A}, j \leq N_n\}$ is bounded from below by a term in the form,

$$\alpha_n^{-Ck_n}.$$

The prior probability of $S_{n,2} = \{|\varepsilon_j - \varepsilon'_j| \leq \alpha_n^{-2\beta-1}\varepsilon_j(1-\varepsilon_j), j \leq N_n\}$ is bounded from below by a term in the form,

$$\prod_{j=1}^{N_n}[\varepsilon_j(1-\varepsilon_j)]^T \alpha_n^{2N_n(\beta+1)} \geq \alpha_n^{-N_n[2(\beta+1)+T]}.$$

Since $N_n \leq CN_0\sqrt{\alpha_n}(\log \alpha_n)^{3/2}$, for all $\alpha \in [c_1\alpha_n, c_2\alpha_n]$, together with the condition on $\alpha_\alpha$, we obtain that there exists $C_1 > 0$ independent of $N_n$ such that

(4.5) $\quad \pi(B^*(f_0, \tau_n, p)) \geq e^{-N_n C_1 \log n} c \geq e^{-C_1 N_0 \sqrt{\alpha_n}(\log \alpha_n)^{5/2}}.$

Set $\alpha_n = \alpha_0 n^{2/(2\beta+1)}(\log n)^{-5/(2\beta+1)}$, then $\tau_n \geq \tau_0 n^{-\beta/(2\beta+1)}(\log n)^{5\beta/(4\beta+2)} = \varepsilon_n$.

We now determine an upper bound on the entropy on some sieve of the support of the *adaptive prior*. Denote $\alpha_{0n} = e^{-n^{1/(2\beta+1)}(\log n)^{5\beta/(2\beta+1)}}$ and $\alpha_{1n} = \alpha_0 \times n^{2/(2\beta+1)}(\log n)^{5\beta/(2\beta+1)}$, and set

$$\mathcal{F}_{n,a} = \{(P, \alpha); k \leq k'_n, \alpha_{0n} \leq \alpha \leq \alpha_{1n}; \varepsilon_j > \varepsilon_0, \forall j\}$$

with $\alpha_0, c > 0$, $k'_n = k'_1 n^{1/(2\beta+1)}(\log n)^{q_\beta}$ with $q_\beta = 5\beta/(2\beta+1)$ if $L(k) = 1$ and $q_\beta = (3\beta-1)/(2\beta+1)$ if $L(k) = \log k$ in the definition of the prior on $k$, and $\varepsilon_0$ is defined by

$$\varepsilon_0 = \exp\{-an^{1/(2\beta+1)}(\log n)^{5\beta/(2\beta+1)}\}.$$

Since $\pi_\alpha$ is bounded, for some $c > 0$,

$$\pi(\mathcal{F}_{n,a}^c) \leq e^{-cn\varepsilon_n^2}.$$

To bound the entropy on $\mathcal{F}_{n,a}$, we use Lemma C.1 with the following parametrization: Write $a = \alpha/(1-\varepsilon)$, $a' = \alpha'/(1-\varepsilon')$, $b = \alpha/\varepsilon$ and $b' = \alpha'/\varepsilon'$ and consider $\rho > 0$ small enough, then if $|a'-a| \leq \tau_1 < a$ and $|b'-b| \leq \tau_2 < b$,

$$g_{\alpha',\varepsilon'}(x) \leq \frac{x^{a-\tau_1-1}(1-x)^{b-\tau_2-1}}{B(a-\tau_1, b-\tau_2)} \frac{B(a-\tau_1, b-\tau_2)}{B(a', b')},$$

so that
$$\frac{B(a-\tau_1, b-\tau_2)}{B(a', b')} \leq 1 + \tau_n \quad \Rightarrow \quad |g_{\alpha',\varepsilon'} - g_{\alpha,\varepsilon}| \leq \tau_n.$$

Consider first $\alpha < 2\varepsilon \wedge (1-\varepsilon)$. If

(4.6) $\qquad |\varepsilon - \varepsilon'| \leq \rho\tau_n\varepsilon(1-\varepsilon), \qquad |\alpha - \alpha'| \leq \rho\tau_n\alpha,$

then using case (i) of Lemma C.1 and simple algebra, we obtain
$$|g_{\alpha',\varepsilon'} - g_{\alpha,\varepsilon}| \leq 4\rho\tau_n.$$

We now consider the $\alpha, \varepsilon$'s such that $2(1-\varepsilon) < \alpha < 2\varepsilon$. If

(4.7) $\qquad |\varepsilon - \varepsilon'| \leq \rho \frac{\tau_n \varepsilon(1-\varepsilon)}{\log(\alpha/(1-\varepsilon))}, \qquad |\alpha - \alpha'| \leq \frac{\alpha\rho\tau_n}{\log(\alpha/(1-\varepsilon))},$

then using case (ii) of Lemma C.1 and simple algebra, we obtain
$$|g_{\alpha',\varepsilon'} - g_{\alpha,\varepsilon}| \leq 2\rho'\tau_n$$

for some $\rho' > 0$. Last we consider the case where $\alpha > 2\varepsilon \vee (1-\varepsilon)$. If

(4.8) $\qquad |\varepsilon - \varepsilon'| \leq \frac{\rho\tau_n\varepsilon^2(1-\varepsilon)^2}{\log(\alpha/\varepsilon(1-\varepsilon))}, \qquad |\alpha - \alpha'| \leq \frac{\rho\varepsilon(1-\varepsilon)\tau_n}{\alpha\log(\alpha/\varepsilon(1-\varepsilon))},$

then case (iv) of Lemma C.1 implies
$$|g_{\alpha',\varepsilon'} - g_{\alpha,\varepsilon}| \leq 2\rho'\tau_n$$

for some $\rho' > 0$. Therefore, the number of intervals in $\alpha$ needed to cover $(e^{-n^{1/(2\beta+1)}(\log n)^{5\beta/(2\beta+1)}} \leq \alpha \leq \alpha_0 n^{2/(2\beta+1)}(\log n)^{10\beta/(2\beta+1)})$ is bounded by
$$J_1 \leq Cn^D \varepsilon_1^{-1} \leq Cn^D e^{(a+1)n^{1/(2\beta+1)}(\log n)^{5\beta/(2\beta+1)}},$$

where $C, D$ are positive constants. We now consider the entropy associated with the supporting points of $P$. The most restrictive relation is (4.8).

Let $\varepsilon_{n,j} = \varepsilon_0^{1/j}$, $j = 1, \ldots, J$ with
$$J = \frac{an^{1/(2\beta+1)}(\log n)^{5\beta/(2\beta+1)}}{t \log n} = \frac{ak'_n}{k'_1 t},$$

so that $\varepsilon_{n,J} = n^{-t}$. Let $P = \sum_{i=1}^{k} p_i g_{\alpha,\varepsilon_i}$, and $N_{n,j}(P)$ be the number of points in the support of $P$ belonging to $(\varepsilon_{n,j}, \varepsilon_{n,j+1})$.

The number of intervals following relation (4.8) needed to cover $(\varepsilon_{n,j}, \varepsilon_{n,j+1})$ is bounded by
$$J_{n,j} = \frac{[\log(\varepsilon_{n,j+1}) - \log(\varepsilon_{n,j})]n^{D_1}}{\varepsilon_{n,j}}$$

for some positive constant $D_1$ independent of $t$. The number of intervals following relation (4.8) needed to cover $(n^{-t}, 1/2)$ is bounded by $J_{n,J+1} = n^{t+1}(\log n)^q$ for some positive constant $q$. For simplicity's sake we consider $D_1 = D_2$. We index the interval $(n^{-t}, 1/2)$ by $J+1$. Consider a configuration $\sigma$ in the form $N_{n,j}(P) = k_j$, for $j = 1, \ldots, J+1$ where $\sum_j k_j = k \leq k'_n$, and define $\mathcal{F}_{n,a}(\sigma) = \{P \in \mathcal{F}_{n,a}; N_{n,j}(P) = k_j, j = 1, \ldots, J+1\}$. For each configuration, the number of balls needed to cover $\mathcal{F}_{n,a}(\sigma)$ is bounded by $J_n(\sigma) = \prod_{j=1}^{J+1} J_{n,j}^{k_j}$. Moreover, the prior probability of $\mathcal{F}_{n,a}(\sigma)$ is bounded by

$$\pi(\mathcal{F}_{n,a}(\sigma)) \leq \Gamma(k+1) \prod_{j=1}^{J+1} \frac{p_{n,j}^{k_j}}{\Gamma(k_j+1)}, \qquad p_{n,j} \leq c[\varepsilon_{n,j+1}^{T+1} - \varepsilon_{n,j}^{T+1}], j \leq J,$$

for some positive consistent $c > 0$ and $p_{n,J+1} \leq 1$. We obtain, since $T \geq 1$ and $t > 2$

$$\Delta_n = \sum_\sigma \sqrt{\pi(\mathcal{F}_{n,a}(\sigma))} \sqrt{J_n(\sigma)}$$

$$\leq \Gamma(k+1)^{1/2} \sum_\sigma \frac{n^{(t+1)k_{J+1}/2}}{\Gamma(k_{J+1}+1)^{1/2}} \prod_{j=1}^J (Cn^{D_1})^{k_j/2} \frac{\varepsilon_{n,j+1}^{(T+1)k_j/2}}{\varepsilon_{n,j}^{k_j/2} \Gamma(k_j+1)^{1/2}}$$

$$\times [\log(\varepsilon_{n,j+1}) - \log(\varepsilon_{n,j})]^{k_j/2} \left[1 - \frac{\varepsilon_{n,j}^{T+1}}{\varepsilon_{n,j+1}^{T+1}}\right]^{k_j/2}.$$

Since

$$\prod_{j=1}^J \Gamma(k_j+1)^{1/2} \leq \exp(k \log(k+1)) \leq e^{k \log(n)},$$

if $tT > 6$, we have

$$\Delta_n \leq C^k n^{kD_1} \Gamma(k+1)^{1/2} \sum_\sigma \prod_{j=1}^J \frac{\exp\{-ak_j k'_n \log n[Tj-2]/(2k'_1 j(j+1))\}}{\Gamma(k_j+1)^{1/2}}$$

$$\leq C^k n^{kD_1} \Gamma(k+1)^{1/2} \sum_\sigma \prod_{j=1}^J \frac{\exp\{-tTk_j \log n/3\}}{\Gamma(k_j+1)^{1/2}}$$

$$\leq C^k n^{k(D_1+t/2+1/2)} \Gamma(k+1)^{1/2} \exp\left\{-\frac{tTk \log n}{6}\right\} \sum_\sigma \prod_{j=1}^J \frac{1}{\Gamma(k_j+1)}$$

$$\leq C^k n^{k(D_1+t/2+1/2)} \Gamma(k+1)^{-1/2} \exp\left\{-\frac{tTk \log n}{6}\right\} e^{k \log(J)}$$

$$\leq e^{k(D_1 - Tt/6 + t/2 + 1/2) \log n}.$$

Hence by choosing $\tau_n = \tau_0 n^{\beta/(2\beta+1)}(\log n)^{q_\beta+1}$ with $\tau_0$ large enough, the above term multiplied by $e^{-n\tau_n^2}$ goes to 0 with $n$, which achieves the proof of Theorem 2.1.

4.2. *Proof of Theorem 2.2.* The proof for the control of the prior mass of Kullback–Leibler neighborhoods of the true density under the *Dirichlet prior* follows the same line as the proof under the *adaptive prior*. To find a lower bound on $\pi(B^*(f_0, \tau_n, p))$, we construct a subset of $\pi(B^*(f_0, \tau_n, p))$ whose probability under a Dirichlet process is easy to compute. Consider $\alpha \in (c_1\alpha_n, c_2\alpha_n)$ and the discrete distribution $P(\varepsilon) = \sum_{j=0}^{N_n} p_j \delta_{\varepsilon_j}(\varepsilon)$ with $N_n = N_0\sqrt{\alpha}(\log \alpha)^{3/2}$ and $\alpha^{-t_0} = \varepsilon_0 < \varepsilon_1 < \cdots < \varepsilon_{N_n} = 1 - \alpha^{-t_0}$ and such that

$$\mathrm{KL}(f_0, g_{\alpha,P}) \leq C\alpha^{-\beta}, \qquad V_p(f_0, g_{\alpha,P}) \leq C\alpha^{-\beta}.$$

The above computations [leading to (4.4)] imply that there exists $D_1$ such that if $|\varepsilon - \varepsilon'| < \alpha^{-D_1}$, we can replace $g_{\alpha,\varepsilon}$ by $g_{\alpha,\varepsilon'}$ in the expression of $g_{\alpha,P}$ without changing the order of approximation of $f_0$ by $g_{\alpha,P}$. Hence we can assume that the point masses $\varepsilon_j$ of the support of $P$ satisfy $|\varepsilon_j - \varepsilon_{j+1}| \geq \alpha^{-D_1}$, $j = 0, \ldots, N_n$. We can thus construct a partition of $[\varepsilon_0/2, 1 - \varepsilon_0/2]$, namely $U_0, \ldots, U_{N_n}$, with $\varepsilon_j \in U_j$ and $\mathrm{Leb}(U_j) \geq 2^{-1}\alpha_n^{-D_1}$ for all $j = 1, \ldots, N_n$ where Leb denotes the Lebesgue measure. Let $\rho > 0$ and $P_1$ be any probability on $[0, 1]$ satisfying

(4.9)       $|P_1(U_j) - p_j| \leq p_j \alpha^{-\rho} \qquad \forall j = 0, \ldots, N_n.$

Then $P_1[\varepsilon_0/2, (1 - \varepsilon_0/2)] \geq 1 - \alpha^{-\rho}$. Since

$$g_{\alpha,P_1}(x) \geq \tilde{g}_{n,P_1} = \int_{\varepsilon_0/2}^{1-\varepsilon_0/2} g_{\alpha_n,\varepsilon}(x)\, dP_1(\varepsilon)$$

and using (4.1), we obtain

$$\mathrm{KL}(f_0, g_{\alpha,P_1}) \leq C\alpha^{-\beta} + \int f_0(x) \log(g_{\alpha,P}(x)/\tilde{g}_{n,P_1}(x))\, dx.$$

Set $\rho \geq \beta$, then, similarly to before, we obtain inequality (4.2) with $\tilde{g}_{n,P_1}$ instead of $g_{\alpha_n,P'}$. When $x \leq e^{-\beta_0\alpha}$, we use the calculations leading to (4.4) replacing $g_{\alpha,P'}$ with $\tilde{g}_{n,P_1}$, wich finally leads to (4.4) between $g_{\alpha,P}$ and $\tilde{g}_{n,P_1}$. To bound $\int_0^1 f_0(x) |\log(\frac{f_0(x)}{g_{\alpha,P_1}(x)})|^p\, dx$, note first that

$$g_{\alpha,P_1}(x) - \tilde{g}_{n,P_1}(x) \leq P_1[0, \varepsilon_0] C\sqrt{\alpha}(x(1-x))^{-1} \leq C\alpha^{-\rho+1/2}(x(1-x))^{-1}.$$

For the purpose of symmetry, we work on $[0, 1/2]$ and we split $[0, 1/2]$ into $[0, e^{-\beta_0\alpha}]$ $[e^{-\beta_0\alpha}, \varepsilon_0]$ $[\varepsilon_0, 1/2]$. Since

$$g_{\alpha,P}(x) \geq g_{\alpha,f_1} - |g_{\alpha,f_1} - g_{\alpha,P}| \geq \frac{f_0(x)}{4} - \frac{C\alpha^{-H}}{x(1-x)} \qquad \forall H > 0,$$

when $x \in (\varepsilon_0, 1/2)$, we have $g_{\alpha,P}(x) \geq cf_0(x)$, since $f_0(x) \geq C_0 x^{k_0}$ near the origin, for some positive constant $c$. Hence combining the above inequality with (4.2) based on $g_{\alpha,P}$ and $\tilde{g}_{n,P_1}$, we obtain that

$$\frac{C\alpha^{-\rho+1/2}}{\tilde{g}_{n,P_1}(x)x(1-x)} \leq C \frac{\alpha^{-\rho+1/2}}{f_0(x)x(1-x)}$$
$$\leq C\alpha^{-\rho+1/2+(k_0+1)t_0}$$
$$= O(\alpha^{-\beta/p}), \quad \text{if } x \in (\varepsilon_0, 1/2) \rho \geq \beta/p + 1/2 + (k_0+1)t_0.$$

Moreover, (4.2) implies, also, that for all $x \in (e^{-\beta_0\alpha_n}, \alpha_n^{-t_0})$,

$$\tilde{g}_{n,P_1} \geq g_{\alpha,P}(x)/2 \geq (x/\varepsilon_0)^\alpha \frac{C\sqrt{\alpha}F_0(\varepsilon_0)}{x(1-x)},$$

leading to

$$\log\left(1 + \frac{g_{\alpha,P_1}(x) - \tilde{g}_{n,P_1}(x)}{\tilde{g}_{n,P_1}(x)}\right) \leq \log\left(1 + \frac{C\alpha^{-\rho+1/2}\varepsilon_0^{\alpha-k_0-1}}{x^\alpha\sqrt{\alpha}}\right)$$
$$\leq C\alpha|\log(x)| \quad \forall x \in (e^{-\beta_0\alpha}, \alpha^{-t_0}).$$

Also, if $x < e^{-\beta_0\alpha}$, using similar calculations to those used to derive (4.3), we obtain

$$\tilde{g}_{n,P_1} \geq CP_1([\varepsilon_0, 1/2])\frac{\sqrt{\alpha}}{x(1-x)}e^{2\alpha\log(x)}$$

and

$$\log\left(1 + \frac{g_{\alpha,P_1}(x) - \tilde{g}_{n,P_1}(x)}{\tilde{g}_{n,P_1}(x)}\right) \leq C\alpha|\log(x)| \quad \forall x < e^{-\beta_0\alpha}.$$

Finally, we obtain

$$\int_0^1 f_0(x)\left|\log\left(\frac{f_0(x)}{g_{\alpha,P_1}(x)}\right)\right|^p dx \leq O(\alpha^{-\beta}) + \int_0^1 f_0(x)\left|\log\left(\frac{\tilde{g}_{n,P_1}(x)}{g_{\alpha,P_1}(x)}\right)\right|^p dx$$
$$\leq \alpha_n^{-t_0+p}(\log\alpha)^p + O(\alpha^{-\beta}) = O(\alpha^{-\beta}),$$

whenever $t_0 > \beta + p$, which implies $\rho > \beta/p + 1/2 + (\beta + p)(k_0 + 1)$.

Under the *Dirichlet prior*, $(P_1(U_0), P_1(U_1), \ldots, P_1(U_{N_n}))$ follows a Dirichlet $(\nu(U_0), \nu(U_1), \ldots, \nu(U_{N_n}))$ with $U_0$ being the complementary set of $(U_1 \cup \cdots \cup U_{N_n})$. Using the fact that $\nu(U_j) \geq C\alpha_n^{-T_1 D_1}$ for all $j$, we obtain that there exist $D_2, C_2 > 0$ such that

$$\pi(S_n) \geq \exp\{-D_2 N_n \log(\alpha_n)\} \geq e^{-C_2 N_0 \sqrt{\alpha_n}(\log\alpha_n)^{5/2}}.$$

The above inequality can be derived, for instance, from Lemma A.2 of [4]. Setting $\alpha_n = n^{2/(2\beta+1)}(\log n)^{-5/(2\beta+1)}$ implies that $\tau_n \geq \tau_0 n^{-\beta/(2\beta+1)}(\log n)^{5\beta/(4\beta+1)} = \varepsilon_n$.

We now bound the $L_1$ entropy for the *Dirichlet prior*. To do so, we use the approximation of any mixture of Beta densities by a finite mixture, and we bound the entropy of a finite mixture. We cannot use the control of the entropy for the *adaptive prior*, however, since it is based on the prior mass of partitions of sieves $\mathcal{F}_{n,a}(\sigma)$, which is not easily controlled under Dirichlet priors. Let $\varepsilon_0 = \exp\{-a\sqrt{\alpha_n}(\log \alpha_n)^{5/2}\}$, $\alpha_n$ as above, and define

$$\mathcal{F}_n = \{F; F[\varepsilon_0, 1-\varepsilon_0] > 1 - \alpha_n^{-\beta}, n^t \leq \alpha \leq \alpha_n(\log \alpha_n)^5\}.$$

Under a Dirichlet $\nu$-process,

$$\pi((\mathcal{F}_n)^c) \leq \alpha_n^\beta \left[\frac{\nu[0, \varepsilon_0]}{\nu[0,1]} + \frac{\nu[1-\varepsilon_0, 1]}{\nu([0,1])}\right] + \exp\{-b\sqrt{\alpha_n}(\log \alpha_n)^{5/2}\}$$

$$\leq C\alpha_n^\beta \exp\{-a\sqrt{\alpha_n}(\log \alpha_n)^{5/2}\}.$$

For all $F \in \mathcal{F}_n$, define $F_n$, the renormalized restriction of $F$, on $[\varepsilon_0, 1-\varepsilon_0]$. Then

$$\|g_{\alpha_n, F_n} - g_{\alpha_n, F}\|_1 \leq 2\alpha_n^{-\beta}.$$

We can, therefore, assume that $F[\varepsilon_0, 1-\varepsilon_0] = 1$ for all $F \in \mathcal{F}'_n$. Then there exists a discrete probability

$$(4.10) \qquad P(\varepsilon) = \sum_{j=1}^{N_n} p_j \delta_{\varepsilon_j}(\varepsilon), \qquad \varepsilon_j \in (\varepsilon_0, 1-\varepsilon_0) \ \forall j,$$

with $N_n \leq N_0\sqrt{\alpha}(\log \alpha)^{3/2}$ such that (3.16) is satisfied for $F$ for all $H$ (depending on $N_0$), and

$$\int_0^{\varepsilon_0/3} |g_{\alpha, F} - g_{\alpha, P}|(x)\, dx \leq \int_{\varepsilon_0}^{1-\varepsilon_0} [dF(\varepsilon) + dP(\varepsilon)] \left(\int_0^{\varepsilon_0/2} g_{\alpha, \varepsilon}(x)\, dx\right).$$

When $x < \varepsilon_0/3 < \varepsilon/3$, using (A.5) we obtain

$$g_{\alpha, \varepsilon}(x) \leq \frac{C\sqrt{\alpha}}{x(1-x)} \left(\frac{2x}{x+\varepsilon}\right)^{\alpha\varepsilon/(2\varepsilon(1-\varepsilon))} \leq C\sqrt{\alpha}\varepsilon^{-\alpha/(2(1-\varepsilon))}(2x)^{\alpha/(2(1-\varepsilon))-1},$$

which implies that

$$\int_0^{\varepsilon_0/3} g_{\alpha, \varepsilon}(x)\, dx \leq C\alpha^{-1/2}(1-\varepsilon)\varepsilon^{-\alpha/(2(1-\varepsilon))} \left(\frac{\varepsilon_0}{3}\right)^{\alpha/(2(1-\varepsilon))}$$

$$(4.11) \qquad \leq C\alpha^{-1/2}(1-\varepsilon)(3/2)^{-\alpha/(2(1-\varepsilon))}$$

$$= O(\alpha^{-H}) \qquad \forall H > 0.$$

By symmetry, the same bound is obtained for the integral over $(1-\varepsilon_0/2, 1)$. Finally, for all $H > 0$, there exists $N_0 > 0$ and a probability measure $P$ defined by (4.10) with $N_n = N_0\sqrt{\alpha}(\log \alpha)^{3/2}$ such that

$$\|g_{\alpha,F} - g_{\alpha,P}\|_1 \leq \alpha^{-H}.$$

Hence the entropy of $\mathcal{F}_n$ is bounded by the entropy of the set

$$\mathcal{F}'_n = \left\{P = \sum_{j=1}^{k} p_j g_{\alpha,\varepsilon_j}; k \leq N_n; \varepsilon_j \in (\varepsilon_0, 1-\varepsilon_0) \ \forall j; n^t \leq \alpha \leq \alpha_n(\log \alpha_n)^5\right\}.$$

Let $k \leq N_n$ be fixed and $g_{\alpha,P}$ be a Beta mixture with $k$ components. When $|\varepsilon'_j - \varepsilon_j| \leq \delta\alpha^{-\gamma_1-2}\varepsilon_j(1-\varepsilon_j)$ for all $j \leq k$ and $|p_j - p'_j| \leq \alpha^{-\gamma_1-1}$, if $|x - \varepsilon_j| \leq \varepsilon_j(1-\varepsilon_j)M_\alpha$, then Lemma B.1 implies

$$|g_{\alpha,\varepsilon'_j} - g_{\alpha,\varepsilon_j}| \leq g_{\alpha,\varepsilon_j} C \alpha^{-\gamma_1}\sqrt{\log \alpha_n},$$

and if $|x - \varepsilon_j| > \varepsilon_j(1-\varepsilon_j)M_\alpha$, then $|x - \varepsilon'_j| > \varepsilon_j(1-\varepsilon_j)\frac{M_\alpha}{2}$ and the convexity of $x \to K(\varepsilon, x)$ for all $\varepsilon$, together with (3.4), implies

$$|g_{\alpha,\varepsilon'_j} + g_{\alpha,\varepsilon_j}| \leq C\frac{\alpha}{x(1-x)}e^{-M^2 \log \alpha/12}.$$

Combining the above inequality with (4.11) leads to

$$(4.12) \qquad \int_0^1 |g_{\alpha,\varepsilon_j} - g_{\alpha,\varepsilon'_j}|(x)\, dx = O(\alpha^{-\gamma_1})$$

and

$$\int_0^1 |g_{\alpha,P} - g_{\alpha,P'}|(x)\, dx = O(\alpha_n^{-\beta}),$$

by choosing $\gamma_1$ large enough. Similarly, considering $|\alpha - \alpha'| \leq n^{-B}\alpha$, we obtain, using (3.5),

$$|g_{\alpha,\varepsilon}(x) - g_{\alpha',\varepsilon}(x)| \leq Cg_{\alpha,\varepsilon}(x)n^{-B},$$

leading to

$$\int_0^1 |g_{\alpha,P} - g_{\alpha',P'}|(x)\, dx = O(\alpha_n^{-\beta}),$$

by choosing $B$ large enough. The number of balls needed to cover $(n^t, \alpha_n(\log \alpha_n)^5)$ under the above constraint is bounded by $Cn^B \log n$. The number of balls with radii $\delta_1\alpha_n^{-\gamma_1}$ needed to cover the set $\mathcal{S}_k$ is bounded by

$$C^k \alpha_n^{-k\gamma_1}.$$

The number of balls with radii $\varepsilon_j(1-\varepsilon_j)\alpha_n^{-\gamma_1-2}\delta_0$ needed to cover $(\varepsilon_0, 1-\varepsilon_0)$

$$(a\alpha_n^{5/2+\gamma_1}(\log \alpha_n)^{5/2})^k.$$

Finally, the metric entropy is bounded by

$$\mathcal{J}_n(\tau_n) \leq 3k_n\beta\log\alpha_n \leq 3k_1\beta\sqrt{\alpha_n}(\log\alpha_n)^5 \leq Cn\tau_n^2,$$

which achieves the proof of Theorem 2.2.

## APPENDIX A: PROOF OF LEMMA 3.2

Throughout the proof, $C$ denotes a generic constant. Let

$$I_0(x) = g_\alpha(x) - 1 = \int_0^1 g_{\alpha,\varepsilon}(x)\,d\varepsilon - 1.$$

The aim is to approximate $I_0$ with an expansion of terms in the form $Q_j(x)\alpha^{-j/2}$ where $Q_j$ is a polynomial function. The idea is to split the integral into three parts, $I_1, I_2, I_3$ corresponding to $\varepsilon < x - \delta_x$, $\varepsilon > x + \delta_x$ and $|x - \varepsilon| < \delta_x$ where $\delta_x = \delta_0 x(1-x)\sqrt{\log(\alpha)/\alpha}$, for some well-chosen $\delta_0 > 0$. Note that this choice of $\delta_x$ comes from the approximation of the Beta density with a Gaussian with mean $x$ and variance $x^2(1-x)^2/\alpha$. We first prove that the first two parts are very small and the expansion is obtained from the third term. By convexity of $K(\varepsilon, x)$ as a function of $\varepsilon$, $K(\varepsilon, x) \geq K(x-\delta_x, x)$ for all $\varepsilon < x - \delta_x$, and $K(\varepsilon, x) \geq K(x+\delta_x, x)$ for all $\varepsilon > x + \delta_x$. Moreover,

$$K(x-\delta_x, x) = x\left(1 - \frac{\delta_0(1-x)\sqrt{\log(\alpha)}}{\sqrt{\alpha}}\right)\log\left(1 - \frac{\delta_0(1-x)\sqrt{\log(\alpha)}}{\sqrt{\alpha}}\right)$$
$$+ (1-x)\left(1 + x\frac{\delta_0\sqrt{\log(\alpha)}}{\sqrt{\alpha}}\right)\log\left(1 + x\frac{\delta_0\sqrt{\log(\alpha)}}{\sqrt{\alpha}}\right)$$
$$= \frac{\delta_0^2\log(\alpha)x(1-x)}{2\alpha} + O\left(x(1-x)\left(\frac{\log(\alpha)}{\alpha}\right)^{3/2}\right),$$

uniformly in $x$. Using a similar argument on $K(x+\delta_x, x)$, we finally obtain, when $\alpha$ is large enough,

$$(A.1) \qquad K(x-\delta_x, x) \geq \frac{\delta_x^2}{3x(1-x)}, \qquad K(x+\delta_x, x) \geq \frac{\delta_x^2}{3x(1-x)}.$$

Set

$$I_1(x) = \int_0^{x-\delta_x} g_{\alpha,\varepsilon}(x)\,d\varepsilon.$$

First we consider $x \leq 1/2$, then using (3.4) and the fact that if $\alpha$ is large enough, the term in the square brackets in (3.4) with $k=1$ is bounded by 2, uniformly in $\varepsilon$, we obtain that

$$I_1(x) \leq \frac{2\sqrt{\alpha}}{\sqrt{2\pi}x(1-x)} \int_0^{x-\delta_x} e^{-\delta_0^2 x(1-x)\log\alpha/(3\varepsilon(1-\varepsilon))} \, d\varepsilon.$$

Let $\rho = (\delta_0^2 x(1-x)\log\alpha)/6$ then

(A.2) $$I_1(x) \leq \frac{2\sqrt{\alpha}}{\sqrt{2\pi}} e^{-\rho/(x-\delta_x)} \leq C\sqrt{\alpha} e^{-\delta_0^2 \log\alpha/6}.$$

Now we consider $x > 1/2$, for which we use another type of upper bound: we split the interval $(0, x-\delta_x)$ into $(0, x(1-\delta))$ and $(x(1-\delta), x-\delta_x)$ for some well-chosen positive constant $\delta$. For all $\varepsilon < x(1-\delta)$, $K(\varepsilon, x) \geq K(x(1-\delta), x)$. Since $u \log(u)$ goes to zero when $u$ goes to zero, there exists $\delta_1 > 0$ such that for all $x > 1/2$, and all $\delta_1 < \delta < 1$,

$$K(x(1-\delta), x) = x(1-\delta)\log(1-\delta)$$
$$+ (1 - x + \delta x)\log\left(1 + \frac{\delta x}{1-x}\right)$$
$$\geq \delta^2 x \log\left(1 + \frac{\delta x}{1-x}\right).$$

Therefore, using (3.4) and the same bound on the square brackets term in (3.4), as in the case $x \leq 1/2$, we obtain that if $x > 1/2$,

(A.3)
$$\int_0^{x(1-\delta)} g_{\alpha,\varepsilon}(x) \, d\varepsilon$$
$$\leq \frac{\sqrt{\alpha}}{\sqrt{2\pi}x(1-x)} \int_0^{x(1-\delta)} \left(1 + \frac{\delta}{2(1-x)}\right)^{-\alpha\delta^2/(2\varepsilon(1-\varepsilon))} d\varepsilon$$
$$\leq \frac{C\sqrt{\alpha}}{(1-x)}\left(1 + \frac{\delta}{2(1-x)}\right)^{-\alpha\delta^2/2}$$
$$\leq C\alpha^{-H} \quad \forall H > 0.$$

We now study the integral over $(x(1-\delta), x-\delta_x)$. We use the following lower bound on $K(\varepsilon, x)$: a Taylor expansion of $K(\varepsilon, x)$ as a function of $\varepsilon$ around $x$ leads to

$$K(\varepsilon, x) = \varepsilon \log\left(\frac{\varepsilon}{x}\right) + (1-\varepsilon)\log\left(\frac{1-\varepsilon}{1-x}\right)$$
$$= (\varepsilon - x)^2 \int_0^1 \frac{(1-u)}{(x + u(\varepsilon - x))(1 - x - u(\varepsilon - x))} \, du$$

(A.4)
$$\geq \frac{(\varepsilon - x)^2}{2} \int_0^{1/2} \frac{1}{(1 - x + u(x - \varepsilon))} \, du$$
$$= \frac{(x - \varepsilon)}{2}(\log(1 - x/2 - \varepsilon/2) - \log(1 - x)).$$

Let $u = x - \varepsilon$, and note that the function $u \to u/(x-u)(1-x+u)$ is increasing so that when $\alpha$ is large enough, uniformly in $x$,

$$g_{\alpha,\varepsilon}(x) \leq \frac{2\sqrt{\alpha}}{\sqrt{2\pi}x(1-x)}\left(\frac{1-x+u/2}{1-x}\right)^{-\alpha u/(2(x-u)(1-x+u))}$$
$$\leq \frac{2\sqrt{\alpha}}{\sqrt{2\pi}x(1-x)}\left(\frac{1-x+u/2}{1-x}\right)^{-\alpha\delta/(2(1-\delta)(1-x+\delta x))}$$

for all $u \in (\delta_x, \delta x)$. Thus if $\alpha$ large enough and $x > 1/2$,

$$\int_{x(1-\delta)}^{x-\delta_x} g_{\alpha,\varepsilon}(x) \, d\varepsilon$$
$$\leq \frac{2\sqrt{\alpha}}{\sqrt{2\pi}x(1-x)} \int_{\delta_x}^{\delta x} \left(1 + \frac{u}{2(1-x)}\right)^{-\alpha\delta/(2(1-\delta)(1-x+\delta x))}$$
$$\leq \frac{8\sqrt{\alpha}}{\sqrt{2\pi}} \frac{1}{\alpha\delta/(2(1-\delta)(1-x+\delta x)) - 1}\left(1 + \frac{\delta_x}{2(1-x)}\right)^{-\alpha\delta/(2(1-\delta)(1-x+\delta x))}$$
$$\leq \frac{C}{\sqrt{\alpha}} e^{-\delta\delta_0\sqrt{\alpha}\sqrt{\log(\alpha)}/(2(1-\delta))} = o(\alpha^{-H})$$

for any $H > 0$. Finally, the above inequality, together with (A.3) for $x > 1/2$ and with (A.2) for $x \leq 1/2$ implies that

$$I_1(x) = O(\alpha^{-H})$$

for all $H > 0$ by choosing $\delta_0$ large enough. We now consider the integral over $(x + \delta_x, 1)$

$$I_2(x) = \int_{x+\delta_x}^{x(1+\delta)} g_{\alpha,\varepsilon}(x) \, d\varepsilon + \int_{x(1+\delta)}^{1} g_{\alpha,\varepsilon}(x) \, d\varepsilon.$$

First let $x \leq 1/2$, then when $\varepsilon \in (x + \delta_x, x(1+\delta))$ with $\delta$ small enough, we can use (3.6) and

$$\int_{x+\delta_x}^{x(1+\delta)} g_{\alpha,\varepsilon}(x) \, d\varepsilon \leq 2e^{-\delta_0^2 \log \alpha/2}.$$

When $\varepsilon \in (x(1+\delta), 1)$, a Taylor expansion of $K(\varepsilon, x)$ as a function of $\varepsilon$ around $x$ leads to

$$K(\varepsilon, x) = (\varepsilon - x)^2 \int_0^1 \frac{(1-u)}{(x + u(\varepsilon - x))(1 - x - u(\varepsilon - x))} du$$

(A.5)
$$\geq \frac{(\varepsilon - x)^2}{2} \int_0^{1/2} \frac{1}{(x + u(\varepsilon - x))} du$$

$$= \frac{(\varepsilon - x)}{2} (\log((x+\varepsilon)/2) - \log x).$$

Thus letting $u = \varepsilon - x$ and noting that $\varepsilon(1-\varepsilon) \leq x + u$ and that $u/(x+u) \geq \delta/(1+\delta)$ as soon as $u > \delta x$, we obtain

$$\int_{x(1+\delta)}^1 g_{\alpha,\varepsilon}(x) \, d\varepsilon \leq \frac{C\sqrt{\alpha}}{x} \int_{\delta x}^{1-x} \left(\frac{2x}{2x+u}\right)^{\alpha\delta/(2(1+\delta))} du$$

$$\leq 2C\alpha^{-1/2} \left(1 + \frac{\delta}{2}\right)^{-\alpha\delta/(2(1+\delta))+1}.$$

If $x > 1/2$ and $\varepsilon > x + \delta_x$, by symmetry, we obtain the same result as in the case $x \leq 1/2$ and $\varepsilon < x - \delta_x$ changing $x$ into $1 - x$. Finally, choosing $\delta_0$ large enough, we prove that for all $x \in [0, 1]$,

(A.6)
$$I_1(x) + I_2(x) = o(\alpha^{-H})$$

($H$ depending on $\delta_0$). We now study the last term, $I_3(x)$. Using (3.6), and the fact that when $\varepsilon \in (x - \delta_x, x + \delta_x)$,

$$|R(x, \varepsilon)| \leq R' \alpha^{k_2+1} |x - \varepsilon|^{3(k_2+1)} (x(1-x))^{-3(k_2+1)}$$

$$\leq R' \alpha^{-(k_2+1)/2} (\log \alpha)^{3(k_2+1)/2},$$

we obtain, for all $k_2 \geq 1, k_1 \geq 3(k_2 - 1)$, and considering the change of variable $u = \sqrt{\alpha}(x - \varepsilon)/(x(1-x))$,

$$I_3(x) = \int_{x-\delta_x}^{x+\delta_x} g_{\alpha,\varepsilon} \, d\varepsilon - 1$$

$$= \sum_{j=1}^{k_2} \frac{\mu_{3j} C(x)^j}{\alpha^{j/2}} + \sum_{j=1}^{k_2 k_1} \frac{\mu_j B_j(x)}{\alpha^{j/2}} + O(\alpha^{-(k_2+1)/2} (\log \alpha)^{3(k_2+1)/2})$$

$$= \frac{I(x)}{\alpha} + O(\alpha^{-(k_2+1)/2} (\log \alpha)^{3(k_2+1)/2}),$$

choosing $\delta_0$ large enough, and since $\mu_1 = 0$ where the $B_j$'s are polynomial functions of $x$ coming from $Q_{k_1}$ and $C(x)$ and where the remaining term is uniform in $x$. Lemma 3.2 is proved.

## APPENDIX B: LEMMA B.1

LEMMA B.1. *Let $(\delta_n)_n$, $(\beta_n)_n$ and $(\rho_n)_n$ be positive sequences decreasing to 0 and assume that $\alpha_n$ increases to infinity. Let $1 - \delta_n > \varepsilon, \varepsilon' > \delta_n$ and $|\varepsilon - \varepsilon'| \leq \rho_n \varepsilon(1-\varepsilon)/\sqrt{\alpha_n}$, then for all $|x - \varepsilon| \leq M\varepsilon(1-\varepsilon)\sqrt{\log \alpha_n}/\sqrt{\alpha_n}$, if $\rho_n\sqrt{\log \alpha_n}$ goes to 0 as n goes to infinity, for all $k_2, k_3 > 1$,*

$$\left| \frac{g_{\alpha_n,\varepsilon}(x)}{g_{\alpha_n,\varepsilon'}(x)} - 1 \right| \leq C[\rho_n\sqrt{\log \alpha_n} + \alpha_n^{-k_2/2}(\log \alpha_n)^{k_2/2} + \alpha_n^{-k_3}]$$

*for n large enough. Also, for all $x \in (\beta_n, 1-\beta_n)$, if $\alpha_n^{1/2}\rho_n|\log(\beta_n)|\delta_n^{-1} = o(1)$, for n large enough,*

$$\left| \frac{g_{\alpha_n,\varepsilon}(x)}{g_{\alpha_n,\varepsilon'}(x)} - 1 \right| \leq C[\alpha_n^{1/2}\rho_n|\log(\beta_n)|\delta_n^{-1} + \alpha_n^{-k_2/2}(\log \alpha_n)^{-k_2/2} + \alpha_n^{-k_3}].$$

PROOF. First let $|x-\varepsilon| \leq M\varepsilon(1-\varepsilon)\sqrt{\log \alpha_n}/\sqrt{\alpha_n}$, since $|\varepsilon-\varepsilon'| \leq \rho_n\varepsilon(1-\varepsilon)/\sqrt{\alpha_n}$, we have that

$$|x - \varepsilon'| \leq \varepsilon(1-\varepsilon)\alpha_n^{-1/2}[M\sqrt{\log \alpha_n} + \rho_n]$$
$$\leq 2M\varepsilon(1-\varepsilon)\alpha_n^{-1/2}\sqrt{\log \alpha_n}$$

and

$$(x-\varepsilon')^l = (x-\varepsilon)^l + (\varepsilon - \varepsilon')\sum_{i=1}^{l} C_l^i (\varepsilon - \varepsilon')^{i-1}(x-\varepsilon)^{l-i}$$
$$= (x-\varepsilon)^l + O(\alpha_n^{-l/2}\rho_n\varepsilon^l(1-\varepsilon)^l(\log \alpha_n)^{(l-1)/2}).$$

We control $g_{\alpha_n,\varepsilon}/g_{\alpha_n,\varepsilon'}$ using approximation (3.5). Then noting that when $n$ is large enough,

$$\left| 1 + \frac{(x-\varepsilon)}{x(1-x)}\left(C(x) + Q_{k_1}\left(\frac{x-\varepsilon}{x(1-x)}\right)\right) \right| \leq 2$$

and

$$\alpha_n|\varepsilon - \varepsilon'||x - \varepsilon| \leq 2x^2(1-x)^2\rho_n\alpha_n^{1/2}(\log \alpha_n)^{1/2},$$

we obtain that

$$a_n = \left| \frac{\alpha_n(x-\varepsilon)^2}{2x^2(1-x)^2}\left[1 + \frac{(x-\varepsilon)}{x(1-x)}\left(C(x) + Q_{k_1}\left(\frac{x-\varepsilon}{x(1-x)}\right)\right)\right] \right.$$
$$\left. - \frac{\alpha_n(x-\varepsilon')^2}{2x^2(1-x)^2}\left[1 + \frac{(x-\varepsilon')}{x(1-x)}\left(C(x) + Q_{k_1}\left(\frac{x-\varepsilon'}{x(1-x)}\right)\right)\right] \right|$$
$$\leq C[\rho_n^2 + (\log \alpha_n)^{1/2}\rho_n + (\log \alpha_n)\rho_n\alpha_n^{-1/2}]$$

and finally,

$$\left|\frac{g_{\alpha_n,\varepsilon}}{g_{\alpha_n,\varepsilon'}}(x) - 1\right| \leq C\rho_n\sqrt{\log \alpha_n} + O(\alpha_n^{1-k_2/2}\varepsilon^{k_2}(1-\varepsilon)^{k_2}(\log \alpha_n)^{k_2/2} + \alpha_n^{-k_3}).$$

Now let $|x - \varepsilon| > M\varepsilon(1-\varepsilon)\sqrt{\log(\alpha)}/\sqrt{\alpha_n}$ and $x \in (\beta_n, 1-\beta_n)$, we use (3.4) together with the above calculations and the fact that the function $\varepsilon \to \varepsilon \log(\varepsilon)/(1-\varepsilon)$ is bounded on $[0,1]$,

$$\begin{aligned}\frac{g_{\alpha_n,\varepsilon}}{g_{\alpha_n,\varepsilon'}}(x) &= \exp\left\{-\alpha_n\left[\frac{1}{1-\varepsilon}\log\left(\frac{\varepsilon}{x}\right) - \frac{1}{1-\varepsilon'}\log\left(\frac{\varepsilon'}{x}\right)\right.\right.\\ &\qquad\qquad\left.\left. + \frac{1}{\varepsilon}\log\left(\frac{1-\varepsilon}{1-x}\right) - \frac{1}{\varepsilon'}\log\left(\frac{1-\varepsilon'}{1-x}\right)\right]\right\}\\ &\quad \times (1 + O(\rho_n\alpha_n^{-1} + \alpha_n^{-k_3}))\\ &= \exp\left\{-\alpha_n(\varepsilon-\varepsilon')\left[\frac{\tilde\varepsilon}{1-\tilde\varepsilon}\log(\tilde\varepsilon) - \frac{1-\tilde\varepsilon}{\tilde\varepsilon}\log(1-\tilde\varepsilon)\right.\right.\\ &\qquad\qquad\left.\left. - \log(x)\frac{\tilde\varepsilon}{(1-\tilde\varepsilon)} - \log(1-x)\frac{(1-\tilde\varepsilon)}{\tilde\varepsilon}\right]\right\}\\ &\quad \times (1 + O(\rho_n\alpha_n^{-1} + \alpha_n^{-k_3})),\end{aligned}$$

where $\tilde\varepsilon \in (\varepsilon, \varepsilon')$. Hence as soon as $1 - \delta_n > \varepsilon, \varepsilon' > \delta_n$ and $x \in (\beta_n, 1-\beta_n)$,

$$\left|\log(1-x)\frac{(1-\tilde\varepsilon)}{\tilde\varepsilon}\right| \leq |\log(\beta_n)|\delta_n^{-1}, \qquad \left|\log(x)\frac{\tilde\varepsilon}{(1-\tilde\varepsilon)}\right| \leq |\log(\beta_n)|\delta_n^{-1},$$

which implies that if $\alpha_n^{1/2}\rho_n|\log(\beta_n)|\delta_n^{-1}$ is small enough,

$$\left|\frac{g_{\alpha_n,\varepsilon}}{g_{\alpha_n,\varepsilon'}}(x) - 1\right| \leq C\alpha_n^{1/2}\rho_n|\log(\beta_n)|\delta_n^{-1} + O(\rho_n\alpha_n^{-1} + \alpha_n^{-k_3}),$$

which achieves the proof of Lemma B.1. □

## APPENDIX C: LEMMA C.1

The following lemma allows us to control the ratio of constants of Beta densities.

LEMMA C.1. *Let $a, b > 0$ and $0 < \tau_1 < a$, $0 < \tau_2 < b$, let $C, \rho$ denote generic positive constants. Let $\bar\eta = a + b$ and $\bar\tau = \tau_1 + \tau_2$. We then have the following results:*

(i) *If $a, b < 2$,*

$$\log\left(\frac{\Gamma(a-\tau_1)\Gamma(b-\tau_2)}{\Gamma(a+\tau_1)\Gamma(b+\tau_2)}\right) + \log\left(\frac{\Gamma(\bar\eta+\bar\tau)}{\Gamma(\bar\eta-\bar\tau)}\right) \leq \frac{2\tau_1}{a-\tau_1} + \frac{2\tau_2}{b-\tau_2} - 2(\bar\tau)C.$$

(ii) *If $a < 2$, $b > 2$, then $\bar{\eta} > 2$ and*

$$\log\left(\frac{\Gamma(a-\tau_1)\Gamma(b-\tau_2)}{\Gamma(a+\tau_1)\Gamma(b+\tau_2)}\right) + \log\left(\frac{\Gamma(\bar{\eta}+\bar{\tau})}{\Gamma(\bar{\eta}-\bar{\tau})}\right) \leq \frac{2\tau_1}{a-\tau_1} + \bar{\tau}[\log(\bar{\eta}+1) - C].$$

(iii) *If $b < 2$, $a > 2$, then things are symmetrical to the previous case.*

(iv) *If $a, b > 2$, $i = 1, 2$, then*

$$\log\left(\frac{\Gamma(a-\tau_1)\Gamma(b-\tau_2)}{\Gamma(a+\tau_1)\Gamma(b+\tau_2)}\right) + \log\left(\frac{\Gamma(\bar{\eta}+\bar{\tau})}{\Gamma(\bar{\eta}-\bar{\tau})}\right) \leq 2\bar{\tau}\log(\bar{\eta}+1).$$

PROOF. The proof of Lemma C.1 comes from Taylor expansions of $\log(\Gamma(x))$ and from the use of the relation,

$$\psi(x) = -\frac{1}{x} + \psi(x+1),$$

so that when $x$ is small, $|\psi(x)|$ is bounded by $1/x$ plus a constant, and if $x$ is large, $\psi(x)$ is bounded by $\log(x)$ plus a constant. $\square$

**Acknowledgments.** I thank C. P. Robert for pointing out the richness of Beta mixtures, together with W. Kruijer for giving me insightful ideas on how these mixtures could be handled. This work has been partially supported by the ANR-SP Bayes grant. I also thank the referees and the editor for pointing out interesting references and useful comments.

CEREMADE
Université Paris Dauphine
Place du Maréchal deLattre de Tassigny
75016 Paris
France
E-mail: rousseau@ceremade.dauphine.fr